\documentclass[twocolumn]{autart}    
\usepackage{graphicx}          
\usepackage{amsmath,amsfonts}
\usepackage{xcolor}
\usepackage{amssymb}
\usepackage{subcaption}
\usepackage{physics}

\begin{document}
\bibliographystyle{abbrv}  
\allowdisplaybreaks
\begin{frontmatter}
	
\title{Mean-Square Exponential Stabilization of \\Mixed-Autonomy Traffic PDE System} 
\thanks[footnoteinfo]{Corresponding author: Huan Yu (huanyu@ust.hk)}

\author[gz]{Yihuai Zhang},    
\author[gz,hk]{Huan Yu},               
\author[fr]{Jean Auriol}, 
\author[MikeFr]{Mike Pereira}

\address[gz]{The Hong Kong University of Science and Technology (Guangzhou), Thrust of Intelligent Transportation, Guangzhou, China.}  
\address[hk]{The Hong Kong University of Science and Technology, Department of Civil Engineering, Hong Kong, China.} 
\address[fr]{Universit\'{e} Paris-Saclay, CNRS, CentraleSup\'{e}lec, Laboratoire des Signaux et Syst\`{e}mes, 91190, Gif-sur-Yvette, France}             
\address[MikeFr]{Mines Paris-PSL University, Department of Geosciences and Geoengineering, 77300, Fontainebleau, France}        


\begin{abstract}                          
Control of mixed-autonomy traffic where Human-driven Vehicles (HVs) and Autonomous Vehicles (AVs) coexist on the road has gained increasing attention over the recent decades. This paper addresses the boundary stabilization problem for mixed traffic on freeways. The traffic dynamics are described by uncertain coupled hyperbolic partial differential equations (PDEs) with Markov jumping parameters, which aim to address the distinctive driving strategies between AVs and HVs. Considering that the spacing policies of AVs vary in mixed traffic, the stochastic impact area of AVs is governed by a continuous Markov chain. The interactions between HVs and AVs such as overtaking or lane changing are mainly induced by impact areas. Using backstepping design, we develop a full-state feedback boundary control law to stabilize the deterministic system (nominal system). Applying Lyapunov analysis, we demonstrate that the nominal backstepping control law is able to stabilize the traffic system with Markov jumping parameters, provided the nominal parameters are sufficiently close to the stochastic ones on average. The mean-square exponential stability conditions are derived, and the results are validated by numerical simulations.

\end{abstract}
\end{frontmatter}

\baselineskip=.97 \normalbaselineskip 
\section{Introduction}
Stop-and-go traffic oscillations,  as common traffic instabilities on freeways, lead to increased travel time, fuel consumption, and the risk of traffic accidents. Stabilization of stop-and-go traffic is achieved mainly via traffic management infrastructures such as ramp metering and varying speed limits. The rapid development of Autonomous Vehicles (AVs) in recent years has brought an increasing market penetration rate of AVs on road~\cite{calvert2017will}. It was initially envisioned that AVs could operate with reduced car-following distances, potentially leading to increased road capacity. This notion has been explored through both theoretical modeling and field experiments~\cite{kesting2008adaptive}. However, it is noteworthy that current commercial AVs usually keep longer distance than Human-driven Vehicles (HVs). This is primarily due to safety concerns as AVs require a sufficient safety buffer to react effectively in the risky event, such as the preceding vehicle comes to a complete stop~\cite{li2022trade}.

Macroscopic modeling for pure HV traffic can be represented by the Lighthill and Whitham and Richard (LWR) model~\cite{lighthill1955kinematic,richards1956shock}, and the Aw-Rascle-Zhang(ARZ) model~\cite{aw2000resurrection,zhang2002non}. Among these, the ARZ model is widely adopted for its ability to describe the stop-and-go oscillations, compared to the LWR. The ARZ model comprises coupled first-order hyperbolic PDEs, which describe the evolution of traffic density and velocity states. To stabilize traffic at its equilibrium, numerous studies have attempted to address the stabilization problem using boundary backstepping control~\cite{anfinsen2018adaptive,espitia2022traffic,hu2015control,krstic2008boundary,vazquez2011backstepping,yu2022traffic}. In~\cite{yu2019traffic,yu2022traffic} the authors first proposed boundary feedback control laws to reduce traffic oscillations in congested regimes, developing full-state feedback and output feedback control laws to achieve exponential stability of the linearized ARZ model. However, stochastic disturbances are not considered. The authors in~\cite{auriol2020robust} designed robust output feedback regulation to ensure hyperbolic systems' robustness against delays and uncertainties. PI boundary controllers~\cite{zhang2019pi} are also applied to stabilize the traffic system described by PDEs around the given steady state.
In addition, distributed control schemes also play a crucial role in the stabilization of PDE traffic systems. Bekiaris-Liberis~\cite{bekiaris2020pde} proposed a control law to stabilize traffic flow for an ARZ-type model consisting of both Adaptive Cruise Control (ACC)-equipped and manually driven vehicles. The authors in~\cite{qi2022delay} further considered the input delay and developed an input delay-compensating control law to stabilize mixed traffic scenarios. 

For a mixed-autonomy traffic scenario, HVs and AVs coexist on the road. Although it remains an open question for PDE modeling of the mixed-autonomy traffic systems, there have been studies on multi-class vehicles considering differences in vehicle sizes and maneuvers, such as trucks and motorcycles. Mohan~\cite{mohan2017heterogeneous} proposed a two-class traffic PDE model by introducing the concept of area occupancy to capture interactions between different classes. Burkhardt~\cite{burkhardt2021stop} adopted the two-class traffic PDE model and output-feedback boundary control laws using the backstepping method to achieve exponential stability in $L_2$-sense. However, they are dealing with the deterministic condition of the mixed-autonomy traffic system.  Real-world traffic systems can be influenced by various stochastic factors such as diverse driving behaviors among vehicle types~\cite{albaba2019stochastic,jin2020analysis}. These factors result in stochastic systems with random parameters embedded within the deterministic system. We consider the spacing of AVs as stochastic in our study, rendering the system become stochastic. The different spacing settings are independent and we use a Markov process to describe it.

In this paper, we propose a mean-square exponential control method to stabilize the oscillations of the proposed stochastic mixed autonomy traffic system. 
The stability of stochastic linear hyperbolic PDEs has been widely investigated~\cite{amin2011exponential,bolzern2006almost,do2012continuous,lamare2015switching,prieur2014stability}. 
Prieur~\cite{prieur2014stability} modeled the abrupt changes of boundary conditions as a piecewise constant function and derived sufficient conditions for the exponential stability of the switching system. 
In~\cite{wang2012stochastically}, the authors examined the robust stochastically exponential stability and stabilization of uncertain linear first-order hyperbolic PDEs with Markov jumping parameters.  
Auriol~\cite{auriol2023mean} demonstrated the mean-square exponential stability of coupled hyperbolic systems with Markovian parameters. With further applications in traffic flow, Zhang~\cite{zhang2017stochastic} studied traffic flow control of Markov jump hyperbolic systems. 
Unlike the results in \cite{zhang2017stochastic}, this paper focuses on the stochastic mixed-autonomy traffic system, adopts the backstepping control, and proposes a novel Lyapunov candidate to get the exponential mean-square stability results.
The stability of mixed-autonomy traffic containing CAV platoons and HVs with two switching parameters was analyzed in~\cite{yu2020stability}. 
In this paper, we choose the spacing of AVs as stochastic due to the communication loss and control delays in AVs. These factors make the spacing of AVs could not always stay at their reference values.
We then consider a nominal boundary controller designed with the backstepping method that would stabilize the system if the stochastic parameter was constant, equal to a reference value. The reference value refers to the best situations of AVs, meaning there is no communication loss, control delays, or other factors that could make AVs away from the reference value. Provided the nominal parameters are sufficiently close to the stochastic ones on average, this control law will stabilize the stochastic mixed-autonomy traffic system. 

The main results of this paper are that we propose a mixed-autonomy traffic PDE system with the stochastic disturbance of AVs governed by a Markov process. To the best of the authors' knowledge, this is the first model that describes the macroscopic spatial-temporal dynamics of HVs and AVs in traffic, using an extended ARZ model with Markov-jumping parameters. Then a backstepping boundary control law is designed to achieve mean-square stabilization of the traffic system. We further develop Lyapunov analysis on the stochastic traffic model to obtain the desired robust stability results. The theoretical contribution mainly lies in proposing a novel stochastic Lyapunov candidate and obtaining robust stabilization with the nominal backstepping controller. The significance of this paper also extends to practical applications since it paves the way for modeling and control of stochastic mixed-autonomy traffic scenarios.

The paper is organized as follows: In Section 2, the stochastic factors in mixed-autonomy traffic systems are considered, and thus, the stochastic mixed-autonomy traffic model is developed. The well-posedness of the stochastic system under the control law is also proved. Section 3 introduces the backstepping control law and conducts the Lyapunov analysis for the nominal system to prove the stability of the closed-loop system. Section 4 provides a Lyapunov analysis of the stochastic mixed- autonomy traffic system and derives conditions for the mean-square exponential stability of the system. Numerical simulation results validating the theoretical derivations are presented in Section 5. Finally, Section 6 concludes the paper.

\section{Stochastic mixed-autonomy traffic model}
In this section, we present the extended stochastic ARZ model that describes the mixed-autonomy traffic scenarios. The traditional ARZ model only describes one class of vehicle's velocity and density evolution. When the road comprises multiple types of vehicles, the ARZ model can not accurately describe the mixed-traffic scenario and the interaction between the different types of vehicles. To solve this problem, an extended AR model was proposed in~\cite{burkhardt2021stop} to describe the behavior of two types of vehicles.  
\subsection{General model with stochastic parameters}
In this paper, we assume that the road only consists of two types of vehicles (HVs and AVs). The following equations describe the dynamics of the system
{\begin{align}
    \partial_t \rho_1+\partial_x\left(\rho_1 v_1\right) &=0, \\
    \partial_t\left(v_1-V_{e,1}\left(AO\right)\right)&+v_1 \partial_x\left(v_1 -V_{e,1} \left(AO\right)\right)\notag\\
    &=\frac{V_{e, 1}\left(AO\right)-v_1}{\iota_1}, \\
    \partial_t \rho_2+\partial_x\left(\rho_2 v_2\right) &=0, \\
    \partial_t\left(v_2-V_{e, 2}\left(AO\right)\right)&+v_2 \partial_x \left(v_2-V_{e, 2}\left(AO\right)\right)\notag\\
    &=\frac{V_{e, 2}\left(AO\right)-v_2}{\iota_2},
\end{align}}
with the boundary conditions 
\begin{align}
    \rho_1(0, t) & =\rho_1^*, \label{boundary1}\\
    \rho_2(0, t) & =\rho_2^*, \label{boundary2}\\
    \rho_1(0, t) v_1(0, t)+\rho_2(0, t) v_2(0, t) & =q_1^*+q_2^*, \label{boundary3}\\
    \rho_1(L, t) v_1(L, t)+\rho_2(L, t) v_2(L, t) & =q_1^*+q_2^* + U(t),\label{boundary4}
\end{align}
where the spatial and time {$(x,t)$ variables belong to  $\{ [0,L] \times \mathbb{R}^+\}$. The index $1$ corresponds to the first class of vehicles (HVs) and the index $2$ corresponds to the second class (AVs)}. The variables ${\rho}_1(x,t)$ and  ${\rho}_2(x,t)$ respectively represent the traffic density on the road, while the variables ${v}_1(x,t)$, ${v}_2(x,t)$ respectively correspond to the traffic velocity. {The parameters $\iota_1$, $\iota_2$ are the relaxation time of AVs and HVs for the drivers to adapt their velocity to the desired velocity.}The function $U(t)$ is the control input that has to be designed. Finally, the variables $\rho_1^*>0$, $\rho_2^*>0$ respectively are the equilibrium density, while $q_1^*>0$, $q_2^*>0$ are the equilibrium flow rate. 

\begin{figure}
    \centering
    \includegraphics[width = 0.45 \textwidth]{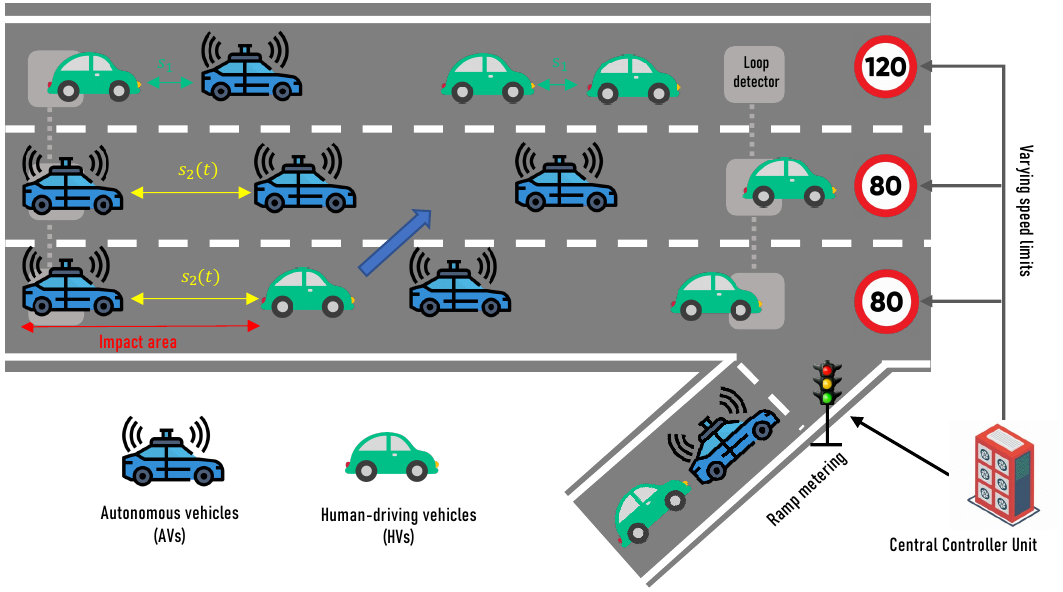}
    \caption{The schematic diagram of the mixed-autonomy traffic control scenario.}
    \label{mixedtraffic}
\end{figure}

The boundary conditions are assumed for an averaged time period of the considered track section, in which the penetration rate of HVs and AVs stay unchanged during the time period and traffic conditions do not drastically change. Thus we assume in \eqref{boundary1}-\eqref{boundary3} the density and flow rate of incoming traffic are at the equilibrium values at $x=0$. The traffic is assumed to be lightly congested and thus the inflow and outflow rates remain the same without control in \eqref{boundary4}. In addition, the control input represents the flow rate perturbation around the steady state and can be implemented by the ramp metering. The ramp inflow actuated by the ramp metering is obtained by  $q_{\text{ramp}}(t) = q_1^\star + q_2^\star + U(t) - q_{\text{main}}(t)$ in which real-time mainline flow $q_{\text{main}}(t)$ can be measured by loop detectors, as shown in Fig. \ref{mixedtraffic}.

Compared with HVs, AVs tend to maintain larger spacing as they adopt more conservative driving strategies \cite{li2022trade}, primarily in response to the safety challenges posed by the current state of autonomous driving technology, which is still evolving. The larger spacing will induce a creeping effect in which HVs will overtake AVs, especially in congested traffic conditions. We set the spacing of AVs as stochastic to reflect its uncertainty. The delays and communication loss among AVs occur stochastically, and they are independent, leading to stochastic driving strategies. Therefore, we choose a continuous Markov chain to describe the stochasticity of AVs. In addition, AVs could have many spacing settings given different driving scenarios. The Markov process with a finite state space guarantees the well-posedness of the model and fits the finite modes of spacing settings. The continuous-time Markov chain allows the stochastic variables to change values at any time not at a set of predefined time steps.

We denote $a_1$ and $a_2$ as impact areas of HVs and AVs:
\begin{align}
    a_1 &= d \times (l + s_1),\\
    a_2 &= d \times (l + s_2(t)),
\end{align}
where $d$ is the vehicle width, $l$ is the vehicle length. The parameter $s_1$ denotes the spacing of HVs, while the function $s_2(t)$ is the spacing of AVs, which will be modeled using a continuous-time Markov process~\cite{zhang2017stochastic}. More precisely, the function $s_2(t)$ can only take a finite number of different values that belong to the set  $\mathfrak{S}=\left\{s^1_2, \ldots, s^r_{2}\right\}$, where $r$ is a positive integer. We consider that the elements of the set $\mathfrak{S}$  are ordered and that there exist known lower bounds and upper bounds $\underline{s}$ and $ \bar{s}$ such that $\underline{s} \leq s^1_2<\cdots<s^{r}_2 \leq \bar{s}$. 
The transition probabilities $P_{ij}(t_1,t_2)$ denote the probability to switch from mode $s_2^i$ at time $t_1$ to mode $s_2^j$ at time $t_2$ ($(i,j) \in \{1, \dots, r\}^2,~0 \leq t_1 \leq t_2$). They satisfy ~$P_{ij}: \mathbb{R}^2 \rightarrow [0,1]$ with $\sum_{j=1}^{r} P_{ij}(t_1,t_2)=1$. Moreover, for $\varrho <t$, we have the Kolmogorov forward equations~\cite{hoyland2009system,kolmanovsky2001mean}:
\begin{align}
    &\partial_t P_{i j}(\varrho, t)=-\mathfrak{c}_j(t) P_{i j}(\varrho, t)+\sum_{k=1}^{r} P_{i k}(\varrho, t) \tau_{k j}(t),\notag\\ 
    & P_{i i}(\varrho, \varrho)=1, \text { and } P_{i j}(\varrho, \varrho)=0 \text { for } i \neq j,
    \label{koleq}
\end{align}
 where the $\tau_{ij}$ and $\mathfrak{c}_j=\sum_{k=1}^{r} \tau_{j k}$ are non-negative-valued functions such that $\tau_{ii}(t)=0$. The functions $\tau_{ij}$ are upper bounded by a constant $\tau^\star$. Finally, we assume that the realizations of $s_2(t)$ are right-continuous~\cite{rausand2003system,ross2014introduction}.
 Based on the definition of impact areas, the interaction of the two classes of vehicles also exists on the road. Hence we introduce area occupancy $AO(\rho_1,\rho_2,t)$~\cite{burkhardt2021stop,mohan2017heterogeneous} to describe the mixed-density of mixed-autonomy traffic:
\begin{align}
    AO(\rho_1,\rho_2, t) = \frac{a_1 \rho_1 + a_2(t) \rho_2}{W},
\end{align}
where $W$ is the road width. Then, we can define the stochastic fundamental diagram combining the area occupancy to describe the relationship between velocity and density of mixed-autonomy traffic using modified Greenshield's fundamental diagram:
\begin{align}
    V_{e,1}(\rho_1,\rho_2,t) = V_1 \left( 1 - \left(\frac{AO}{\overline{AO}_1}\right)^{\gamma_1}\right),\\
    V_{e,2}(\rho_1,\rho_2,t) = V_2 \left( 1 - \left(\frac{AO}{\overline{AO}_2}\right)
    ^{\gamma_2}\right),
\end{align}
where $V_1$, $V_2$ are the free-flow velocity of HVs and AVs, $\overline{AO}_1$, $\overline{AO}_2$ are the maximum area occupancy of them, $\gamma_1$, $\gamma_2$ denote the traffic pressure exponents. 
The fundamental diagrams denote the relationship between the velocity and density of HVs and AVs. Using the stochastic fundamental diagram, we can get the stochastic equilibrium velocity:
\begin{align}
    v_1^*(t) = V_{e,1}(\rho_1^*,\rho_2^*,t),\quad v_2^*(t) = V_{e,2}(\rho_1^*,\rho_2^*,t). \label{Stoc_eq_speed}
\end{align}
Thus, the stochastic equilibrium flow can be calculated by
\begin{align}
    q_1^*(t)= v_1^*(t) \rho_1^*,\quad q_2^*(t) = v_2^*(t) \rho_2^*.
\end{align}
We could also choose $s_1$ as stochastic variable due to the different driving styles of drivers. Once both spacing settings are stochastic, the stochastic process would be a concatenated one. The model can be modified to ensure both spacing settings are stochastic. The paper would become more technical and heavy in notations.

\subsection{Linearization}
Next, we linearize the stochastic traffic flow system around its equilibrium at time $t$ in mode $i$ $\in$ $\{0,\dots,r\}$, denoted as $v_1^*(t)=v_{1i}^*$, $q_{1}^*(t)=q_{1i}^*$.
We define the small deviation from the equilibrium point $\Tilde{\rho}_1 = \rho_1(x,t) - \rho_{1}^*$, $\Tilde{v}_1 = v_1(x,t) - v_{1i}^*$, $\Tilde{\rho}_2 = \rho_2(x,t) - \rho_{2}^*$, $\Tilde{v}_2 = v_2(x,t) - v_{2i}^*$, and define the augmented state $\mathbf{z}(x,t) = \begin{bmatrix}
    \Tilde{\rho}_1(x,t)  &  \Tilde{v}_1(x,t) & \Tilde{\rho}_2(x,t) & \Tilde{v}_2(x,t)
\end{bmatrix}^\mathsf{T}$. Adjusting the computations done in~\cite{burkhardt2021stop}, we obtain the following linearized system
\begin{align}
    \mathbf{z}_t(x,t) + \mathbf{J}_\lambda^i \mathbf{z}_x(x,t) = \mathbf{J}_i\mathbf{z}(x,t),\label{eq_rho_v}
\end{align}
with the boundary conditions
\begin{align}
    \left[ \begin{array}{cccc}
        1 & 0 & 0 & 0 \\
        0 & 0 & 1 & 0\\
        v_{1i}^* & \rho_1^* & v_{2i}^* &\rho_2^* \\
    \end{array} \right] \mathbf{z}(0,t) &= 0,\\
    \left[ \begin{array}{cccc}
        v_{1i}^* & \rho_1^* &v_{2i}^* & \rho_2^* \\
    \end{array} \right] \mathbf{z}(L,t) &= U(t),\label{eq_rho_v_BC}
\end{align}
where
\begin{align*}
    &\mathbf{J}_\lambda^i =\left[\begin{array}{cc}
    v_{1i}^* & \rho_1^* \\
    0 & v_{1i}^*-\beta_{11}^i \rho_1^* \\
    0 & 0 \\
    \beta_{21}^i\left(v_{2i}^*-v_{1i}^*\right) & -\beta_{21}^i \rho_1^*\\
    \end{array}\right. \\
    &\left.\begin{array}{cc}
    0 & 0 \\
    \beta_{12}^i\left(v_{1i}^*-v_{2i}^*\right) & -\beta_{12}^i \rho_2^* \\
    v_{2i}^* & \rho_2^* \\
    0 & v_{2i}^*-\beta_{22}^i \rho_2^*\\
    \end{array} \right], \\
    &\mathbf{J}_i=\left[\begin{array}{cccc}
    0 & 0 & 0 & 0 \\
    -\frac{1}{\iota_1} \beta_{11}^i & -\frac{1}{\iota_1} & -\frac{1}{\iota_1} \beta_{12}^i & 0 \\
    0 & 0 & 0 & 0 \\
    -\frac{1}{\iota_2} \beta_{21}^i & 0 & -\frac{1}{\iota_2} \beta_{22}^i & -\frac{1}{\iota_2}
    \end{array}\right],
\end{align*}
$\beta_{\mathfrak{m}\mathfrak{n}}=\left. - \frac{\partial V_{e,\mathfrak{m}}^i\left(A O\left(\rho_1, \rho_2\right)\right)}{\partial \rho_\mathfrak{n}}\right|_{\rho_1^*, \rho_2^*}, 1 \leq \mathfrak{m},\mathfrak{n} \leq 2.$

\subsection{Analysis of the system~\eqref{eq_rho_v}-\eqref{eq_rho_v_BC} in mode $i$}
\begin{figure}
    \centering
    \includegraphics[width = 0.45\textwidth]{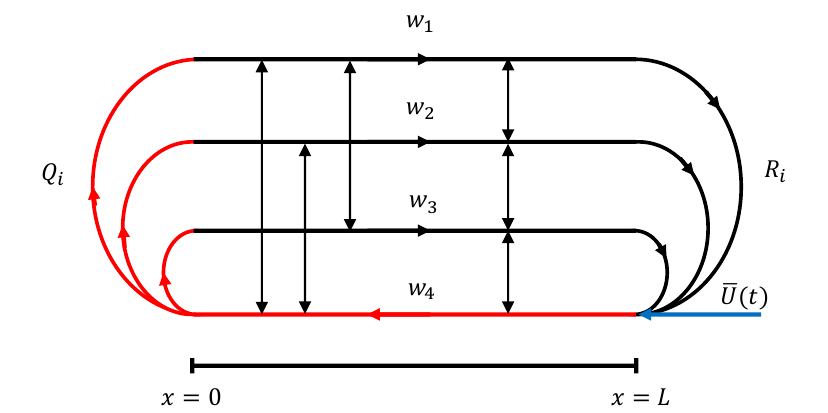}
    \caption{The schematic diagram of system \eqref{clpsys1} - \eqref{clpsys4}: The blue arrow denotes the location where the control input $\Bar{U}(t)$ goes into the system. The boundary couplings between $w_1$, $w_2$ ,$w_3$ and $w_4$ are described by stochastic $Q_i$ and $R_i$.}
    \label{hyperbolicsys1}
\end{figure}

Let us consider for the moment that the system stays in the mode $i \in \{1, \dots r\}$, i.e., that $s_2(t)=s_2^i$, for all $t>0$. Therefore, the system~\eqref{eq_rho_v}-\eqref{eq_rho_v_BC} is not stochastic anymore, and we can apply the change of coordinates introduced in~\cite{burkhardt2021stop} to rewrite it in the Riemann coordinates, with a diagonal transport matrix. Such a representation simplifies the open-loop analysis and will be later required to design our backstepping controller or to derive the proposed Lyapunov analysis.
Let us first define the characteristic velocities $(\lambda^i_1, \lambda^i_2, \lambda^i_3, \lambda^i_4)$ as the eigenvalues of the matrix $\mathbf{J}_{\lambda i}$ corresponding to mode $i$, given by
\begin{align*}
\lambda^i_1 =& v_{1i}^*,\\
\lambda_2^i =& v_{2i}^*,\\
\lambda_3^i =& \frac{v_{1i}^*+v_{2i}^*-\mathcal{V}_{1i}\rho_{1}^*- \mathcal{V}_{2i}\rho_{2}^* + \Delta_i}{2},\\
\lambda_4^i =& \frac{v_{1i}^* +v_{2i}^* -\mathcal{V}_{1i}\rho_{1}^* -\mathcal{V}_{2i}\rho_{2}^* - \Delta_i}{2},\\
\mathcal{V}_{1i}\left(\rho_1^*, \rho_2^*\right) =&\left.-\frac{\partial V_{e,1}^i\left(A O\left(\rho_1, \rho_2\right)\right)}{\partial \rho_1}\right|_{\rho_1^*, \rho_2^*},\\
\mathcal{V}_{2i}\left(\rho_1^*, \rho_2^*\right) =&\left.-\frac{\partial V_{e,2}^i\left(A O\left(\rho_1, \rho_2\right)\right)}{\partial \rho_2}\right|_{\rho_1^*, \rho_2^*},
\end{align*}
where 
\begin{align*}
    \Delta_i\left(\rho_1^*, \rho_2^*\right)=&(\left(\mathcal{V}_{2i} \rho_2^*-\mathcal{V}_{1i} \rho_1^*+v_{1i}^*-v_{2i}^*\right)^2 \nonumber\\
    &+4 (\beta_{12} \beta_{21} \rho_1^* \rho_2^*))^{1/2}.
\end{align*}
It was shown that the eigenvalues satisfy~\cite{zhang2006hyperbolicity} 
\begin{align}
    \lambda_{4}^i \leq \min\{ \lambda_{1}^i,\lambda_{2}^i \} \leq \lambda_{3}^i\leq\max\{\lambda_{1}^i,\lambda_{2}^i\}.
\end{align}
Since $\lambda_{1}^i=v_{1}^i>0$, $\lambda_{2}^i = v_2^i>0$, we have $\lambda_3>0$. Define the matrix $\mathbf{V}_i=\{ \Hat{v}_{k\ell}\}_{1 \leq k, \ell \leq 4}$ as the change of basis matrix such that the coefficient matrix can be diagonalized as $\mathbf{V}^{-1}_i \mathbf{J}_{\lambda}^i \mathbf{V}_i$ = $\text{Diag}\{\lambda_1^i,\lambda_2^i,\lambda_3^i,\lambda_4^i\}$, with  positive eigenvalues in ascending order.  
We choose to re-arrange the eigenvalues of the velocity matrix in ascending order to simplify the backstepping analysis we will use and fit the formalism used in related backstepping papers~\cite{hu2015control}.
We also define the source term matrix as $\Hat{\mathbf{J}}_i = \mathbf{V}^{-1}_{i}{\mathbf{J}}_i \mathbf{V}_i = \{ \Hat{J}_{kl}^i \}_{1 \leq k,\ell \leq 4}$.
Let us define the transformation matrix $\mathbf{T}_i$ as
\begin{align*}
   \mathbf{T}_i = \begin{bmatrix}
        \\
       \mathbf{T}_i^{+}\\
       \\ \hline 
       \mathbf{T}_i^{-}
   \end{bmatrix}= \begin{bmatrix}
       0 & \mathrm{e}^{-\frac{\Hat{J}^i_{22}}{v_{2i}^*} x} & 0 & 0 \\
0 & 0 & \mathrm{e}^{-\frac{\Hat{J}_{33}^i}{\lambda_{3i}} x} & 0 \\
\mathrm{e}^{-\frac{\Hat{J}_{11}^i}{v_{1i}^*} x} & 0 & 0 & 0 \\
0 & 0 & 0 & \mathrm{e}^{-\frac{\Hat{J}_{44}^i}{\lambda_{4i}} x}
\end{bmatrix} \mathbf{V}_i^{-1},
\end{align*}
where $\mathbf{T}_i^{+} \in \mathbb{R}^{3\times 4}$ and $\mathbf{T}_i^{-} \in \mathbb{R}^{1\times 4}$.
We can now define the change of coordinates 
\begin{align}
\left[\begin{array}{c}
w_1,
w_2,
w_3,
w_4
\end{array}\right]^\mathsf{T}=\mathbf{T}_i\mathbf{z}. \label{transfomation}
\end{align}
{Using the notations $\mathbf{w} = [\mathbf{w}^{+}, \mathbf{w}^{-}]$, where $\mathbf{w}^{+} = \mathbf{T}_i^+ \mathbf{z}= [w_1, w_2, w_3]^\mathsf{T}$, $\mathbf{w}^{-} =  \mathbf{T}_i^- \mathbf{z}= w_4$.} This new state verifies the following set of PDEs
{\begin{align}
        \mathbf{w}^+_t(x,t) +\Lambda^{+}_i \mathbf{w}^+_x(x,t) =& \Sigma^{++}_i(x)\mathbf{w}^+ (x,t)\nonumber\\
        &+\Sigma_i^{+-}(x) \mathbf{w}^-(x, t), \label{clpsys1}\\
        \mathbf{w}^-_t(x, t)-\Lambda^{-}_i \mathbf{w}^-_x(x, t) =& \Sigma^{-+}_i(x) \mathbf{w}^+(x,t),\label{clpsys2}
\end{align}}
with the boundary conditions 
{\begin{align}
    \mathbf{w}^+(0,t)  &= {Q}_i \mathbf{w}^-(0, t), \label{clpsys3} \\
    \mathbf{w}^-(L, t) &= {R}_i\mathbf{w}^+(L, t)+\bar{U}(t), \label{clpsys4}
\end{align}}
where the different matrices are defined as
\begin{align*}
\Lambda^{+}_i =&\left[\begin{array}{ccc}
\lambda_{2}^i & 0 & 0 \\
0 & \lambda_{3}^i & 0 \\
0 & 0 & \lambda_{1}^i
\end{array}\right], \Lambda^{-}_i =-\lambda_{4}^i,\\
\Sigma^{++}_i(x) =& \left[\begin{array}{ccc}
0 & \bar{J}_{12}^i(x) & \bar{J}_{13}^i(x) \\
\bar{J}_{21}^i(x) & 0 & \bar{J}_{23}^i(x) \\
\bar{J}_{31}^i(x) & \bar{J}_{32}^i(x) & 0
\end{array}\right],\\
\Sigma^{+-}_i(x) =& \left[\begin{array}{ccc}
\bar{J}_{14}^i(x) & \bar{J}_{24}^i(x) & \bar{J}_{34}^i(x)
\end{array}\right]^{\mathsf{T}},\\
\Sigma^{-+}_i(x) =&\left[\begin{array}{lll}
\bar{J}_{41}^i(x) & \bar{J}_{42}^i(x) & \bar{J}_{43}^i(x)
\end{array}\right],\\
{Q}_i = &-\left[\begin{array}{lll}
    0 & 1 & 0 \\
    0 & 0 & 1 \\
    1 & 0 & 0
    \end{array}\right]\left[\begin{array}{ccc}
         \Hat{v}_{11}^i&  \Hat{v}_{12}^i  & \Hat{v}_{13}^i\\
         \Hat{v}_{41}^i & \Hat{v}_{42}^i & \Hat{v}_{43}^i  \\
         \kappa_{1i} & \kappa_{2i} & \kappa_{3i}
    \end{array}\right]^{-1} \left[ \begin{array}{ccc}
         \Hat{v}_{14}^i  \\
          \Hat{v}_{44}^i \\
         \kappa_{4i}
    \end{array}\right],\\
{R}_i = &-\left[ \begin{array}{cc}
     \frac{\kappa_{2i}}{\kappa_{4i}} \mathrm{e}^{\left(\frac{\Hat{J}_{11}^i}{\lambda_{1i}}-\frac{\Hat{J}_{44}^i}{\lambda_{4i}}\right) L} & \frac{\kappa_{3i}}{\kappa_{4i}} \mathrm{e}^{\left(\frac{\Hat{J}_{22}^i}{\lambda_{2i}}-\frac{\Hat{J}_{44}^i}{\lambda_{4i}}\right) L} 
\end{array} \right.\\
&\left.\begin{array}{c}
    \frac{\kappa_{1i}}{\kappa_{4i}} \mathrm{e}^{\left(\frac{\Hat{J}_{33}^i}{\lambda_{3i}}-\frac{\Hat{J}_{44}^i}{\lambda_{4i}}\right) L}
\end{array}\right],\\
\kappa_{\jmath i}=&v_{1i}^* \Hat{v}_{1 \jmath}^i+\rho_1^* \Hat{v}_{2 \jmath}^i+v_{2i}^* \Hat{v}_{3 \jmath}^i+\rho_2^* \Hat{v}_{4 \jmath}^i, \jmath=1,2,3,4,\\
\bar{U}(t)=&\mathrm{e}^{-\frac{\Hat{J}_{44}^i}{\lambda_{4i}} L} \frac{1}{\kappa_{4i}} U(t),
\end{align*}
where $\Bar{J}^i(x)$ is the corresponding transformation $\mathbf{T}_i$ of $\Hat{J}^i(x)$. The detailed calculations of $\Bar{J}^i(x)$ can be found in \cite{burkhardt2021stop}. 
The properties of system \eqref{eq_rho_v}-\eqref{eq_rho_v_BC} are the same as the ones of system \eqref{clpsys1}-\eqref{clpsys4} due to the invertibility of the transformation that maps system \eqref{eq_rho_v}-\eqref{eq_rho_v_BC} to system \eqref{clpsys1}-\eqref{clpsys4}. The schematic diagram of the system \eqref{clpsys1}-\eqref{clpsys4} is shown in Fig. \ref{hyperbolicsys1}.
The system \eqref{clpsys1}-\eqref{clpsys4} can be rewritten in compact form:
\begin{align}
    \partial_t \mathbf{w}(x, t)+\Lambda_{i} \partial_x \mathbf{w}(x, t)= \Theta_{i}(x) \mathbf{w}(x,t), \label{eq_w_Riemann}
\end{align}
{where $\Lambda_i = \text{Diag}\{\lambda_{2}^i,\lambda_{3}^i,\lambda_{1}^i,\lambda_{4}^i\}$,} with the boundary condition:
\begin{align}
    \left[\begin{array}{l}
\mathbf{w}^{+}(0,t) \\
\mathbf{w}^{-}(L, t)
\end{array}\right]=G_{i}\left[\begin{array}{l}
\mathbf{w}^{+}(L,t) \\
\mathbf{w}^{-}(0, t)
\end{array}\right] + \left[\begin{array}{c}
     0  \\
    \Bar{U}(t)
\end{array}\right],
\end{align}
where the coefficient matrix $\Theta_i$, $G_i$ are:
\begin{align*}
    \Theta_{i} =  \left[\begin{array}{cc}
    \Sigma^{++}_i(x) & \Sigma^{+-}_i(x) \\
    \Sigma^{-+}_i(x) & 0
\end{array}\right], G_{i} = \left[ \begin{array}{cc}
        0 & Q_i \\
        R_i & 0
    \end{array} \right]. 
\end{align*}
The state $\mathbf{w}$ and the original state $\mathbf{z}$ have equivalent $L^2$ norms. Therefore, there exist two constants $m_w>0$ and $M_w>0$ such that
 \begin{align}
    m_w||\mathbf{z}||_{L^2}^2 &\leq ||\mathbf{w}||_{L^2}^2\leq M_w||\mathbf{z}||_{L^2}^2.
\end{align}
Since we have $r$ possible modes for the stochastic parameter $s_2(t)$, it implies that we have $r$ possible configurations for the matrices $G_i, \Lambda_i$ and $\Theta_i$.

\subsection{Open-loop analysis,  well-posedness, and control objective}
The traffic system can be divided into a free-flow regime and a congested regime based on the propagation direction of traffic waves. More details of modeling and analysis can be found in~\cite{burkhardt2021stop,yu2019traffic}.
\begin{itemize}
    \item Free-flow regime: all the eigenvalues $\lambda_{1i}>0$, $\lambda_{2i}>0$, $\lambda_{3i}>0$, $\lambda_{4i}>0$. The traffic oscillations are transporting downstream at corresponding speed $\lambda_{1i}$, $\lambda_{2i}$, $\lambda_{3i}$, $\lambda_{4i}$. No congestion is generated in this condition.
    \item Congested regime: $\lambda_{1i}>0$, $\lambda_{2i}>0$, $\lambda_{3i}>0$, $\lambda_{4i}<0$.
    In the congested regime, the disturbance of traffic flow propagates upstream so that all the vehicles will be influenced, and the road becomes congested.
\end{itemize}
   \begin{figure}[t]
        \centering
        \includegraphics[width = 0.45 \textwidth]{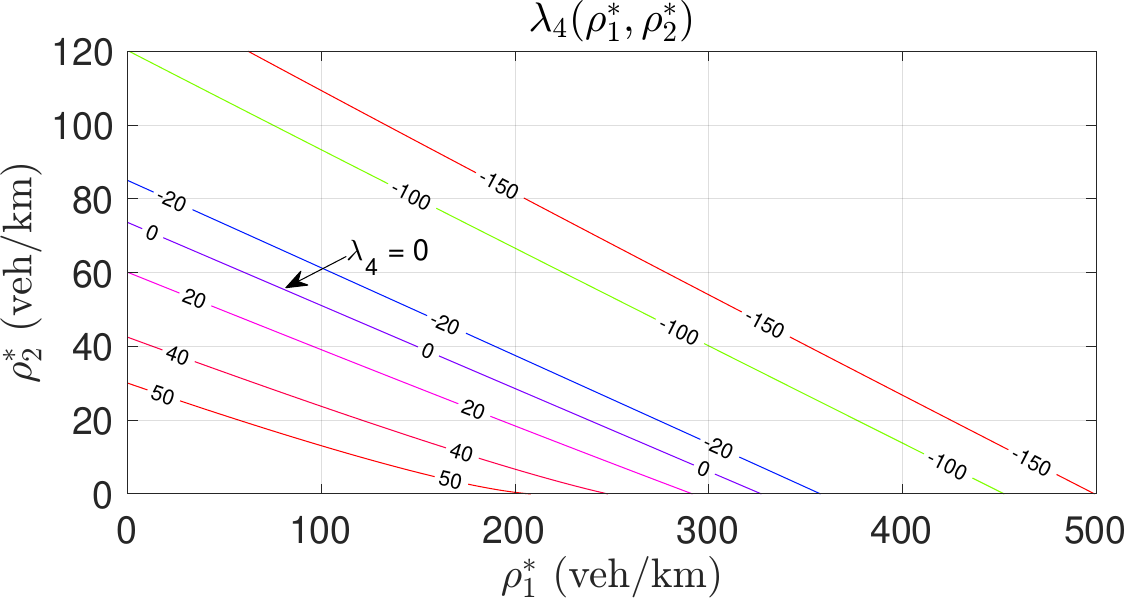}
        \caption{Contour plot of $\lambda_4$}
        \label{lam4to0}
    \end{figure}
In this paper, we are dealing with the traffic in the congested regime, i.e., we have $\lambda_{4i} < 0$ in which both AV and HV have relatively large density values, as shown in Fig. \ref{lam4to0}. Therefore, for all the modes induced by the stochastic parameters, we need to have $\lambda_{4i}<0$ and $\lambda_{1i}>0$, $\lambda_{2i}>0$, $\lambda_{3i}>0$.

In this paper, we will design a controller with the following structure
\begin{align}
    {U}(t)=-F (\mathbf{T}_0^+(\mathbf{z}))(L,t) +\int_0^L K(\xi) (\mathbf{T}_0(\mathbf{z}))(L,t) d\xi. \label{control_law_nominal}
\end{align}
where $F$ is a constant matrix, $K$ is a piecewise continuous function, and where the transformation $\mathbf{T}$ is defined by equation~\eqref{transfomation}. The index $0$ stands for a reference value $s_2^0$ that does not necessarily belong to the set $\mathfrak{S}$. It corresponds to a reference and arbitrary choice for $s_2$ that will be used to design the nominal stabilizing control law.
We now state the well-posedness of the stochastic system with the control law~\eqref{control_law_nominal}
\begin{lem}
For any initial conditions of the Markov system $\mathbf{z}(x,t) \in L^2[0,L]$ and any initial states $s_2(0)$ for the stochastic parameter, the system \eqref{eq_rho_v}-\eqref{eq_rho_v_BC} with the nominal control law \eqref{control_law_nominal} has a unique solution such that for any $t$,
\begin{align}
    \mathbb{E}\{ ||\mathbf{z}(x,t)||_{L^2(0,L)} \} < \infty,
\end{align}\label{unique}
where the $\mathbb{E}\{\cdot\}$ denotes the expectation.
\end{lem}
\vspace{-0.3cm}
\begin{pf}
Almost every sample path of the stochastic process $s_2(t),t\geq 0$ is a right-continuous function with a finite number of jumps in any finite time interval. We can find a sequence $\{ t_k: k =0,1,\dots \}$ of stopping times such that $t_0 = 0, \lim_{t\to \infty} t_k = \infty$, and  $s_2(t) = s_2(t_k)$ on $t_k \leq t < t_{k+1}$. 
Let us consider that $s_2(t_0)=s_{2 i_0} \in \mathfrak{S}$. On the time interval $[t_0,t_1]$, after performing the transformation $\mathbf{T}_{i_0}$, we have
\begin{align}
    \partial_t \mathbf{w}(x, t)+\Lambda_{i_0} \partial_x \mathbf{w}(x, t)= \Theta_{i_0} \mathbf{w}(x,t),
\end{align}
with the boundary condition:
\begin{align}
    \left[\begin{array}{l}
\mathbf{w}^{+}(0,t) \\
\mathbf{w}^{-}(L, t)
\end{array}\right]=G_{i_0}\left[\begin{array}{l}
\mathbf{w}^{+}(L,t) \\
\mathbf{w}^{-}(0, t)
\end{array}\right]+ \left[\begin{array}{c}
     0  \\
    \Bar{U}(t_0)
\end{array}\right].
\end{align}
The initial condition $\mathbf{w}(x,0) \in L^2(0,L)$. {This Cauchy problem has one and only one solution $\mathbf{w}(x, t) \in L^2(0,L)$ using ~\cite[Theorem A.4 ]{bastin2016stability} and \cite[Appendix. A]{coron2021boundary}, such that there exists a constant $C_1$ such that for all $t\in [t_0,t_1]$,}
\begin{align}
    ||\mathbf{w}(x, t)||_{L^2(0,L)} \leq C_1 ||\mathbf{w}(x,0)||_{L^2(0,L)}.
\end{align}
This yields,
\begin{align}
    \mathbb{E}\{||\mathbf{w}(x, t_1)||_{L^2(0,L)}\} <  \infty.
\end{align}
Since $\mathbf{w}(x,t) =\mathbf{T}_{i_0} \mathbf{z}(x,t)$, we obtain
\begin{align}
    \mathbb{E}\{||\mathbf{z}(x, t_1)||_{L^2(0,L)}\} <  \infty.
\end{align}
Let us now consider the next Markov event and the time interval $[t_1, t_2]$. We have $s_2(t)=s_{2 i_1} \in \mathfrak{S}$ and the real system can be written as:
\begin{align}
    \partial_t \mathbf{w}(x, t)+\Lambda_{i_1} \partial_x \mathbf{w}(x, t)= \Theta_{i_1} \mathbf{w}(x,t),
\end{align}
with boundary conditions:
\begin{align}
    \left[\begin{array}{c}
\mathbf{w}^{+}(0,t) \\
\mathbf{w}^{-}(L, t)
\end{array}\right]=G_{i_1}\left[\begin{array}{c}
\mathbf{w}^{+}(L,t) \\
\mathbf{w}^{-}(0, t)
\end{array}\right]+ \left[\begin{array}{c}
     0  \\
    \Bar{U}(t_1)
\end{array}\right].
\end{align}
The initial condition at this time is $\mathbf{w}(x,t_1) \in  L^2(0,L)$. {Using~\cite[Theorem A.4 ]{bastin2016stability} and \cite[Appendix. A]{coron2021boundary} again, we can get the existence of a constant $C_2$ such that the unique solution $\mathbf{w}(x,t)$ satisfies for all $t\in [t_1,t_2]$}
\begin{align}
    ||\mathbf{w}(x, t)||_{L^2(0,L)} \leq C_2||\mathbf{w}(x,t_1)||_{L^2(0,L)},
\end{align}
which implies
\begin{align}
    \mathbb{E}\{||\mathbf{w}(x, t_2)||_{L^2(0,L)}\} <  \infty.
\end{align}
Since  $\mathbf{w}(x,t) = \mathbf{T}_{i_1} \mathbf{z}(x,t)$, we directly obtain
\begin{align}
    \mathbb{E}\{||\mathbf{z}(x, t_2)||_{L^2(0,L)}\} <  \infty.
\end{align}
Iterating the process on the whole time domain, we can get that the original system  $\mathbf{z}(x,t)$ with Markov jump has a unique solution on $[0, \infty)$ that satisfies
\begin{align}
    \mathbb{E}\{ ||\mathbf{z}(x,t)||_{L^2(0,L)} \} < \infty.
\end{align}
This concludes the proof of Lemma~\ref{unique}. $\hfill\blacksquare$
\end{pf}
The linearized open-loop system may be unstable due to the negative transport coefficients and the in-domain coupling terms. When the stochastic parameter is constant, a stabilizing control law has been designed in~\cite{burkhardt2021stop} using the backstepping approach. In this paper, we will consider the design of such a stabilizing controller for a nominal value of $s_2(t)$. We will then show that this control law guarantees the mean-square exponential stability, provided that the difference between the stochastic parameter and this nominal value is small on average.

\section{Nominal backstepping boundary controller and mean-square exponential stabilization }
In this section, we design a nominal controller to stabilize the system if the stochastic parameter $s_2(t)$ always equals a constant reference value $s_2^0$, for all $t \geq 0$. This reference value $s_2^0$ does not necessarily belong to the set $\mathfrak{S}$. It corresponds to a reference and arbitrary choice for $s_2$ to design a nominal stabilizing control law. {In practice, $s_2^0$ denotes the ideal situation for AVs without any communication loss and control delays. AVs drive perfectly on the road in the nominal mode.}
The real values of $s_2$ can then be seen as stochastic disturbances acting on this reference value. We still consider that $\underline{s} < s_2^0 < \bar s$. The corresponding stochastic parameters (that depend on $s_2$ due to different transformations) are also denoted with the index $0$, e.g., $R_0$, $Q_0$. The control design relies on the backstepping approach. From now on, we consider that we are in the nominal mode $s_2(t)=s_2^0$.

\subsection{Backstepping transformation}
Our objective is to simplify the structure of the system~\eqref{eq_w_Riemann} by removing the in-domain coupling terms in the equation. More precisely, let us consider the following backstepping transformation $\mathcal{K}_0$ 
\small
\begin{align}
    \mathcal{K}_0\mathbf{w}= \begin{pmatrix} \mathbf{w}^{+} \\
    \mathbf{w}^- -\int_0^x \mathbf{K}_0(x, \xi)\mathbf{w}^{+}(\xi,t)
    +N_0(x, \xi) \mathbf{w}^-(\xi, t) d \xi \end{pmatrix}\label{back}
\end{align}
\normalsize
where the kernels $\mathbf{K}_0(x,\xi) \in \mathbb{R}^{1\times 3}$ and $N_0(x,\xi)\in\mathbb{R}^1$ are piecewise continuous functions defined on the triangular domain $\mathcal{T}=\{0 \leq \xi \leq x \leq L\}$. We have
\begin{align}
    \mathbf{K}_0(x,\xi) = \left[ \begin{array}{ccc}
        k^0_1(x,\xi) & k^0_2(x,\xi) & k^0_3(x,\xi)
    \end{array}\right].
\end{align}
The different kernels are governed by the following PDEs
\vspace{-0.8 cm}
\small
\begin{align}
&-\Lambda^{-}_0(\mathbf{K}_0)_x(x,\xi) +(\mathbf{K}_0)_\xi(x,\xi)\Lambda^{+}_0=-\mathbf{K}_0(x,\xi)\Sigma^{++}_0(\xi) \notag\\
&-\Sigma^{-+}_0(\xi) N_0(x, \xi), \label{ker1} \\
&\Lambda^{-}_0 (N_0)_{x}(x, \xi)+\Lambda^{-}_0 (N_0)_{\xi}(x, \xi) =\mathbf{K}_0(x,\xi)\Sigma^{+-}_0(\xi), \label{ker2}
\end{align}
\normalsize
with the boundary conditions: 
\begin{align}
\left(-\Lambda^{-}_0 \mathbf{I}_3-\Lambda^{+}_0\right) \mathbf{K}_0^\mathsf{T}(x,x)&=\Sigma^{-+^{\mathsf{T}}}_0(x) ,\label{ker3}\\
-\Lambda^{-}_0 N_0(x, 0)+  \mathbf{K}_0(x,0)\Lambda^{+}_0{Q}_0 &=0, \label{ker4}
\end{align}
where $\mathbf{I}_3$ is a $3\times 3$ identity matrix. The well-posedness of the kernel equations can be proved by adjusting the results from~\cite[Theorem 3.3]{hu2015control}. The solutions of the kernel equations can be expressed by integration along the characteristics. Applying the method of successive approximations, we can then prove the existence and uniqueness of the solution to the kernel equations~\eqref{ker1}-\eqref{ker4}.

Applying the backstepping transformation, we can define the state $\vartheta$ of the target system as $$\vartheta=(\alpha_1,\alpha_2,\alpha_3,\beta)=\mathcal{K}_0 \mathbf{w}.$$ {Using the notation $\delta = [\alpha_1, \alpha_2, \alpha_3]^\mathsf{T}$, the target system equations are given by:}
\begin{align}
&\delta_t(x,t)+\Lambda^{+}_0 \delta_x(x,t) =\Sigma^{++}_0(x)\delta(x,t) +\Sigma^{+-}_0(x) \beta(x,t) \nonumber \\
&+\int_0^x\mathbf{C}_0^+(x,\xi)\delta(\xi,t) d\xi + \int_0^x\mathbf{C}_0^-(x,\xi)\beta(\xi,t)d\xi,
\label{tar1}\\
&\beta_t\left(x, t\right)-\Lambda^{-}_0 \beta_x(x, t)=0, \label{tar2}
\end{align}
with the boundary conditions:
\begin{align}
&\delta(0,t) ={Q}_0 \beta(0, t), \label{tar3}\\
&\beta(L, t)  = R_0 \delta(L,t) + \Bar{U}(t)\notag\\
& - \int_0^L \mathbf{K}_0(L, \xi)\mathbf{w}^{+}(\xi,t)+N_0(L, \xi) \mathbf{w}^-(\xi, t)) d \xi. \label{tar4}
\end{align}
where the coefficients $\mathbf{C}_0^+(x,\xi) \in \mathbb{R}^{3\times 3}$ and $\mathbf{C}_0^-(x,\xi) \in \mathbb{R}^{3\times 1}$ are bounded functions defined on the triangular domain $\mathcal{T}$. {Their expressions are~\cite{burkhardt2021stop}:}
\begin{align*}
        \mathbf{C}_0^+(x,\xi) &= \Sigma^{+-}_0(x)\mathbf{K}_0(x,\xi) + \int_{\xi}^x  \mathbf{C}_0^-(x,s)\mathbf{K}_0(s,\xi) ds\\
        \mathbf{C}_0^-(x,\xi) &=  \Sigma^{+-}_0(x)N_{0}(x,\xi) + \int_\xi^x \mathbf{C}_0^-(x,s)N_0(s,\xi) ds 
\end{align*}
The transformation $\mathcal{K}_0$ is a Volterra type, therefore boundedly invertible~\cite{yoshida1960lectures}. Consequently, the states $\mathbf{w}$ and $\vartheta$ have equivalent $L^2$ norms, i.e. there exist two constants $m_\vartheta>0$ and $M_\vartheta>0$ such that
 \begin{align}
    m_\vartheta||\mathbf{w}||_{L^2}^2 &\leq ||\vartheta||_{L^2}^2\leq M_\vartheta||\mathbf{w}||_{L^2}^2.
\end{align}

\subsection{Nominal control law and Lyapunov functional}
From the nominal target system~\eqref{tar1}-\eqref{tar4}, we can easily design a stabilizing control law as~\cite{auriol2016minimum,redaud2024domain}:
\begin{align}
    \bar{U}(t)=&-{R}_0 \mathbf{w}^+(L,t) +\int_0^L\left(\mathbf{K}_0(L, \xi)\mathbf{w}^+(\xi,t) \right. \notag\\
    & \left. +N_0(L, \xi) \mathbf{w}^- (\xi, t)\right) d \xi. \label{control law}
\end{align}
To analyze the stability properties of the target system~\eqref{tar1}-\eqref{tar4}, we consider the Lyapunov functional $V_0(t)$ defined by 
\begin{align}
     V_0(t) = \int_0^L \vartheta^\mathsf{T}(x,t) D_0(x) \vartheta(x,t)dx, \label{lyapunov_0}
\end{align}
where 
\begin{align}
     D_0(x) = \text{Diag}\left\{\frac{\mathrm{e}^{-\frac{\nu}{\lambda_{2}^0}x}}{\lambda_{2}^0}, \frac{\mathrm{e}^{-\frac{\nu}{\lambda_{3}^0}x}}{\lambda_{3}^0}, \frac{\mathrm{e}^{-\frac{\nu}{\lambda_{1}^0}x}}{\lambda_{1}^0}, a \frac{\mathrm{e}^{\frac{\nu}{\Lambda^-_0} x}}{\Lambda^-_0} \right\}.
    \label{def_D0}
\end{align}
This Lyapunov functional is equivalent to the $L^2$ norm, that is there exist two constant  $k_1 > 0$ and $k_2>0$ such that
\begin{align}
    k_1||\vartheta||_{L^2}^2 \leq V_0(t) \leq k_2||\vartheta||_{L^2}^2. \label{eq_norm_V0}
\end{align}
It can also be expressed in terms of the original state as 
\begin{align}
V_0(t)=\int_0^L&\mathcal{K}_0(\mathbf{T}_0(\mathbf{z}(x,t)))^\mathsf{T} D_0(x)\mathcal{K}_0(\mathbf{T}_0(\mathbf{z}(x,t)))dx.
\end{align}
Taking the time derivative of $V_0(t)$ and integrating by parts, we get:
\begin{align}
    \Dot{V}_0(t) \leq& -\nu V_0(t) + \int_0^L 2  \delta^\mathsf{T}(x,t) D_\alpha^0 \nonumber\\
    &\left(\Sigma_{0}^{++}(x)\delta(x,t) + \Sigma_{0}^{+-}(x) \mathbf{T}_0^-\mathbf{z}(x,t) \right)dx \nonumber \\
    &+(q_{10}^2+q_{20}^2+q_{30}^2 - a) \beta^2(0,t)\nonumber \\
    \leq& - \eta V_0(t) +(q_{10}^2+q_{20}^2+q_{30}^2 - a) \beta^2(0,t),
\end{align}
where 
\begin{align}
    \eta =& \nu - \frac{2}{||\underline{\Lambda^+}|| k_1} (\max_{x\in [0,L]} ||\Sigma^{++}_0(x)||  \nonumber \\
    &+ (1 + \frac{1}{m_\vartheta}) \max_{x\in [0,L]} ||\Sigma^{+-}_0||(x)),\\
    D_\alpha^0=& \text{Diag}\left\{\frac{\mathrm{e}^{-\frac{\nu}{\lambda_{1}^0}x}}{\lambda_{1}^0}, \frac{\mathrm{e}^{-\frac{\nu}{\lambda_{2}^0}x}}{\lambda_{2}^0}, \frac{\mathrm{e}^{-\frac{\nu}{\lambda_{3}^0}x}}{\lambda_{3}^0} \right\}.
\end{align}
We choose $a>0$ and $\nu>0$ such that
\begin{align}
     & q_{10}^2+q_{20}^2+q_{30}^2 - a \leq 0,~\eta>0.
\end{align}
where $q_{10}$, $q_{20}$, $q_{30}$ are the elements of $Q_0$. Consequently, we obtain $\dot{V}_0(t)\leq -\eta V_0(t)$, which implies the $L^2$-exponential stability of the system.
\subsection{Mean-square exponential stabilization}
The main goal of this paper is to prove that the control law (\ref{control law}) can still stabilize the stochastic system \eqref{eq_rho_v}-\eqref{eq_rho_v_BC}, provided the nominal parameter $s_2^0$ is sufficiently close to the stochastic ones on average. More precisely, we have the following robust stabilization result. 
\begin{thm}\label{th_mean_square_stab}
{ There exists a constant $\epsilon^\star>0$, such that for all time $t \geq 0$,}
\begin{align}
  \mathbb{E}\left(\left|s_2^0-s_2(t)\right|\right) \leq \epsilon^\star, \label{eq_ineg_stocha}
\end{align}
then the closed-loop system \eqref{eq_rho_v}-\eqref{eq_rho_v_BC} with the control law (\ref{control law}) is mean-square exponentially stable, namely, there exist $\varsigma,\zeta>0$ such that:
\begin{align}
    \mathbb{E}_{[0,(\mathbf{z}(x,0),s(0))]}(\norm{\mathbf{z}(x,t)}_{L_2}^2) \leq \varsigma \mathrm{e}^{-\zeta t} \norm{\mathbf{z}(x,0)}_{L_2}^2,
\end{align}\label{main_thm}
\end{thm}
This theorem will be proved in the next section.
\begin{rem}
While proving this theorem, we will give a possible value for the bound $\epsilon^\star$ (that always exists). However, this value is likely to be of very little practical use due to the conservativeness of the Lyapunov analysis. Nevertheless, thanks to this expression, one will be able to observe that this positive constant depends on the feedback gain and that this dependence is likely to be considerable. In this context, our result is qualitative: we only guarantee the existence of robustness margins. 
For any constant $\epsilon_1$ that verifies
$\epsilon^\star>\epsilon_1 > 0$ such that $\mathbb{E}\left(\left|s_2^0-s_2(t)\right|\right) \leq \epsilon_1$, it would be interesting to establish the links between the convergence rate and $\epsilon_1$, but this would require a serious quantitative analysis, which is outside the scope of this paper. 
We believe that a smaller $\epsilon_1$ will lead to a faster convergence rate for the stochastic system as the stochastic parameter would be closer on average to the nominal ones. The best situation would correspond to the absence of stochastic terms in the traffic system such that $\epsilon_1$ is 0 and the system corresponds to the nominal system. 
\end{rem}

\section{Lyapunov Analysis}
In this section, we consider the stochastic system with the nominal controller~\eqref{control law}. The objective is to prove Theorem~\ref{th_mean_square_stab}. The proof will rely on a Lyapunov analysis. More precisely, we will consider the following stochastic Lyapunov functional candidate
\begin{align}
    V(t)=\int_0^L&\mathcal{K}_0(\mathbf{T}(t,\mathbf{z}(x,t)))^\mathsf{T} D_j(x)\mathcal{K}_0(\mathbf{T}(t,\mathbf{z}(x,t)))dx, \label{Lyap_functional}
\end{align}
where $(\mathbf{T}(t,\mathbf{z}(x,t))$=$\mathbf{T}_j(\mathbf{z}(x,t))$ if $s_2(t)=s_2^j$. The diagonal matrix $D_j$ is defined by 
\begin{align}
     D_j(x) = \text{Diag}\left\{\frac{\mathrm{e}^{-\frac{\nu}{\lambda_{2}^j}x}}{\lambda_{2}^j}, \frac{\mathrm{e}^{-\frac{\nu}{\lambda_{3}^j}x}}{\lambda_{3}^j}, \frac{\mathrm{e}^{-\frac{\nu}{\lambda_{1}^j}x}}{\lambda_{1}^j}, a \frac{\mathrm{e}^{\frac{\nu}{\Lambda^-_j} x}}{\Lambda^-_j} \right\}.
    \label{def_Dj}
\end{align}
We consider that the parameters $\nu$ and $a$ introduced in the definition of $D_j$ can still be tuned. In the nominal case $s_2(t)=s_2^0$, the Lyapunov functional $V(t)$ corresponds to $V_0(t)$. It is noted that inequality~\eqref{eq_norm_V0} still holds for $V(t)$ (even if the constants $k_1$ and $k_2$ may change). Unlike what has been done in~\cite{auriol2023mean}, this Lyapunov functional is defined in terms of the original state $\mathbf{z}$ instead of the target system state $\vartheta$. This choice is made to avoid taking the derivative of a stochastic signal when using the transformations $\mathbf{T}_i$.

\subsection{Target system in mode $s_2^j$}
In this section, we consider that~$s_2(t)=s_2^j$. We can define the state $\vartheta=(\alpha_1,\alpha_2,\alpha_3,\beta)=\mathcal{K}_0(\mathbf{T}_i(\mathbf{z}))$. Our objective is first to obtain the equations verified by the state $\vartheta$ that appears in the Lyapunov functional~\eqref{Lyap_functional}. Using the same notation $\delta = [\alpha_1, \alpha_2, \alpha_3]^\mathsf{T}$, it verifies the following set of equations
{\begin{align}
&\delta_t(x,t)+\Lambda^{+}_j \delta_x(x,t) =\Sigma^{++}_j(x)\mathbf{T}_j^+\mathbf{z}(x,t) \nonumber\\
&+\Sigma^{+-}_j(x) \mathbf{T}_j^-\mathbf{z}(x,t), \label{tar1_sto}\\
&\beta_t\left(x, t\right)-\Lambda^{-}_j \beta_x(x, t) =\mathbf{f}_{1j}(x) \mathbf{T}^+_j \mathbf{z}(x,t)
+\mathbf{f}_{2j}(x) \beta(0, t) \notag \\
&+\int_0^x \mathbf{f}_{3j}(x,\xi)\mathbf{T}^+_j\mathbf{z}(\xi,t) d \xi 
+\int_0^x \mathbf{f}_{4j}(x,\xi) \mathbf{T}^-_j\mathbf{z}(\xi,t) d \xi, \label{tar2_sto}
\end{align}}
with the boundary conditions:
{\begin{align}
\delta(0,t)  =&{Q}_j \beta(0, t), \label{tar3_sto}\\
\beta(L, t)  =&\left({R}_j-{R}_0\right)\delta(L,t), \label{tar4_sto}
\end{align}}
where the functions are defined by:
\begin{align}
    \mathbf{f}_{1j}(x) =& \Sigma^{-+}_j(x)+\Lambda^-_j \mathbf{K}_0(x, x)+\mathbf{K}_0(x, x) \Lambda^{+}_j,\\
    \mathbf{f}_{2j}(x) =& -\mathbf{K}_0(x, 0) \Lambda^{+}_j {Q}_j+N_0(x, 0) \Lambda^{-}_j,\\
    \mathbf{f}_{3j}(x,\xi) =& \Lambda^{-}_j (\mathbf{K}_0)_{x}(x, \xi)-(\mathbf{K}_0)_{\xi}(x, \xi) \Lambda^{+}_j \notag \\
     &- \mathbf{K}_0(x, \xi) \Sigma^{++}_j(\xi)-N_0(x, \xi) \Sigma^{-+}_j(\xi),\\
    \mathbf{f}_{4j}(x,\xi) =& \Lambda^{-}_j (N_0)_{x} (x, \xi)+\Lambda^{-}_j (N_0)_{\xi}(x, \xi) \notag\\
    &-\mathbf{K}_0(x, \xi) \Sigma^{+-}_j(\xi).
\end{align}
All the terms that depend on $\mathbf{z}$ in the target system~\eqref{tar1_sto}-\eqref{tar4_sto} could be expressed in terms of $\vartheta$ using the inverse transformation~$\mathcal{K}_0^{-1}$. However, this would make the computations more complex and is not required for the stability analysis.
It is important to emphasize that all the terms on the right-hand side of equation~\eqref{tar2_sto} become \emph{small} if the stochastic parameters are close enough to the nominal ones. More precisely, we have the following lemma
\begin{lem}
There exists a constant $M_0$, such that for any realization $s_2(t)=s_{2}^j \in \mathfrak{S}$, for any $(x,\xi)\in \mathcal{T}$
    \begin{align}
        ||\mathbf{f}_{\mathfrak{i}j}|| < M_0  \left|s_{2}^0 - s_{2}^j\right|, \quad \mathfrak{i} \in \{1,2,3,4\}.
    \end{align}\label{lem_bound_f}
\end{lem}
\vspace{-1 cm}
\begin{pf}
First, using the boundedness of the kernel $\mathbf{K}_0$ and the continuity of the matrices $\Sigma_j, \Lambda_j$, such that
\begin{align}
    ||\Sigma^{-+}_j(x) - \Sigma^{-+}_0(x)|| &< C \left|s_2^0 - s_2^j\right|,\notag\\
    ||\Lambda^-_j-\Lambda^-_0|| &< C \left|s_2^0 - s_2^j\right|, \notag\\
    ||\Lambda^{+}_j-\Lambda^+_0|| &< C \left|s_2^0 - s_2^j\right|,
\end{align}
where $C$ is a positive constant.
Considering the function $\mathbf{f}_{1j}(x)$. For all $x \in [0,L]$, we have
\begin{align}
    \mathbf{f}_{1j}(x) =& \Sigma^{-+}_j(x)+\Lambda^-_j \mathbf{K}_0(x, x)+\mathbf{K}_0(x, x) \Lambda^{+}_j\notag \\
    =& (\Sigma^{-+}_j(x) - \Sigma^{-+}_0(x)) + (\Lambda^-_j-\Lambda^-_0) \mathbf{K}_0(x,x) \notag\\
    &+ \mathbf{K}_0(x,x)(\Lambda^{+}_j-\Lambda^+_0).
\end{align}
Consequently, we obtain the existence of a constant $K_1>0$ such that
\begin{align}
    ||\mathbf{f}_{1j}|| &\leq K_1\left|s_2^0 - s_2^{j}\right|.
\end{align}
Consider now the function $\mathbf{f}_{3j}(x,\xi)$. We have:
\begin{align}
    \mathbf{f}_{3j}(x,\xi) &= \Lambda^{-}_j (\mathbf{K}_0)_{x}(x, \xi)-(\mathbf{K}_0)_{\xi}(x, \xi)\Lambda^{+}_j \notag\\
    &- \mathbf{K}(x, \xi) \Sigma^{++}_j(\xi) - \Sigma^{-+}_j(\xi)N(x, \xi), \notag \\
    &\hspace{-0.5cm}= (\Lambda^{-}_j - \Lambda^{-}_0)(\mathbf{K}_0)_x(x, \xi) + (\mathbf{K}_0)_{\xi}(x, \xi)(\Lambda^{+}_0-\Lambda^{+}_j) \notag \\
    &+ \mathbf{K}_0(x,\xi)(\Sigma^{++}_0(\xi) -\Sigma^{++}_j(\xi)) \notag \\
    &+ N_0(x,\xi)(\Sigma^{-+}_0(\xi)- \Sigma^{-+}_j(\xi)).
\end{align}
The regularity of the kernels implies the existence of a constant $K_3>0$ such that 
\begin{align}
    ||\mathbf{f}_{3j}|| &\leq K_3 \left|s_2^0 - s_2^{j}\right|.
\end{align}
{The other inequalities for $\mathbf{f}_{2j}(x)$ and $\mathbf{f}_{4j}(x,\xi)$ can also be derived using similar techniques.} This finishes the proof of Lemma \ref{lem_bound_f}. $\hfill\blacksquare$
\end{pf}
\subsection{Derivation of the Lyapunov function}
Let us consider the Lyapunov functional~ $V$ defined in equation~\eqref{Lyap_functional}. Its infinitesimal generator $L$ is
defined as~\cite{ross2014introduction}
\begin{align}
 L V(\mathbf{z},s_2) =&\limsup _{\Delta t \rightarrow 0^{+}} \frac{1}{\Delta t} \times \mathbb{E}(V(\mathbf{z}(t+\Delta t), s_2(t+\Delta t))\notag\\
&-V(\mathbf{z}(t), s_2(t))).
\end{align}
We define $L_j$, the infinitesimal generator of $V$ obtained  by fixing $s_2(t)=s_2^j$ ($j \in \mathfrak{S}$)
\begin{align}
L_j V(\mathbf{z})  =\frac{d V}{d \mathbf{z}}(\vartheta, s_2^j) h_j(\vartheta) +\sum_{\ell \in \mathfrak{R}}\left(V_{\ell}(\mathbf{z})-V_j(\mathbf{z})\right) \tau_{j \ell}(t),
\end{align}
where $V_\ell(\mathbf{z})=V(\mathbf{z},s_2^\ell)$, and where the operator $h_j$ is defined by 
\begin{align}
    h_j(\vartheta)=\left(\begin{array}{l}
-\Lambda^+_{j} \delta_x(x,t)+\Sigma_{j}^{++}(x)\mathbf{T}_j^+\mathbf{z}(x,t)\\
+\Sigma_{j}^{+-}(x)\mathbf{T}_j^-\mathbf{z}(x,t)\\
\Lambda_{j}^{-} \beta_x(x, t)+\mathbf{f}_{1j}(x)\mathbf{T}_j^+\mathbf{z}(x,t) \\
+ \mathbf{f}_{2j}(x) \beta(0, t)+\int_0^x \mathbf{f}_{3j}(x,\xi)\mathbf{T}_j^+\mathbf{z}(\xi,t) d \xi\\
+\int_0^x \mathbf{f}_{4j}(x,\xi)\mathbf{T}_j^-\mathbf{z}(\xi,t) d \xi
\end{array}\right).
\end{align}
To shorten the computations, we denote in the sequel $V(t), LV(t), V_j(t)$ and $L_jV(t)$ instead of (respectively) $V(\mathbf{z}(t),s_2(t))$, $LV(z(t),s_2(t))$, $V(\mathbf{z}(t),s_2^j)$ and $L_j(V(\mathbf{z}(t)))$. From now, we consider that $s_2(t=0)=s_2^{i_0}$ for some ${i_0}\in\mathfrak{S}$. We have the following lemma
\begin{lem}
There exists $\Bar{\eta}>0$, $M_1 > 0$ and $d_1,d_2 > 0$ such that the Lyapunov functional $V(t)$ satisfies
\begin{align}
&\sum_{j=1}^r P_{i j}(0, t) L_j V(t) \leq -V(t)\Big(\Bar{\eta}-d_1 \mathcal{Z}(t) \nonumber\\
&-\left(M_1+d_1 r \tau^{\star}\Big) \mathbb{E}\left(\left|s_2^0-s_2(t)\right|\right)\right)\nonumber \\
&+\sum_{k=1}^3 (d_2 \mathbb{E}(\left|s_2^0 - s_2(t)\right|)-\mathrm{e}^{-\frac{\nu}{\bar \lambda}})\alpha_k^2(L,t),
\end{align}
where the function $\mathcal{Z}(t)$ is defined as:
\begin{align}
\mathcal{Z}(t)=\sum_{j=1}^r\left|s_2^j-s_2^0\right|\left(\partial_t P_{i j}(0, t)+\mathfrak{c}_j P_{i j}(0, t)\right)
\end{align} \label{lem_lyapunov_functional}
\end{lem}
\vspace{-1 cm}
\begin{pf}
We will first compute the first term of $L_j$. Considering that $s_2(t)=s_2^j$, we can define the state $\vartheta=\mathcal{K}_0(\mathbf{T}_i(\mathbf{z}))$, solution of equations~\eqref{tar1_sto}-\eqref{tar4_sto}. The Lyapunov functional rewrites 
\begin{align}
V_j(t)=\int_0^L\vartheta^T(x,t)D_j(x) \vartheta(x,t)dx,
\end{align}
Direct and classical computations give
{\begin{align}
    &\frac{d V_j}{d \mathbf{z}}( \mathbf{z}) h_j( \mathbf{z}) = -\nu V_j(t) + \int_0^L 2 \delta^{\mathsf{T}}(x,t) D^j_\alpha(x) \notag\\
    &\left(\Sigma_{j}^{++}(x)\mathbf{T}^+_j\mathbf{z}(x,t) + \Sigma_{j}^{+-}(x) \mathbf{T}^-_j\mathbf{z}(x,t) \right)dx \notag \\
    &+\int_0^L \frac{2 a}{\Lambda^-_j} \mathrm{e}^{\frac{\nu}{\Lambda^-_j} x} \beta(x, t)\left(\mathbf{f}_{1j}(x)\mathbf{T}^+_j\mathbf{z}(x,t)+\mathbf{f}_{2j}(x) \beta(0, t)\right. \nonumber\\
    &+\left.\int_0^x \mathbf{f}_{3j}(x,\xi)\mathbf{T}^+_j\mathbf{z}(\xi,t) d \xi \right. \nonumber\\
    &\left.+\int_0^x \mathbf{f}_{4j}(x,\xi) \mathbf{T}^-_j\mathbf{z}(\xi,t) d \xi \right)dx \notag \\
    &+ ({q}_{1j}^2+{q}_{2j}^2+{q}_{3j}^2-a)\beta^2(0,t) - D_{\alpha}^j(L)\delta^2(L,t) \notag \\
    &+ a \mathrm{e}^{\frac{\nu}{\Lambda^-_j}}\left(\left({R}_j-{R}_0\right)\delta(L,t)\right)^2,
\end{align}}
where the ${q}_{1j}$,${q}_{2j}$, ${q}_{3j}$ are the elements of ${Q}_j$ and where $D_\alpha^j$ is defined by
\begin{align}
    D_\alpha^j= \text{Diag}\left\{\frac{\mathrm{e}^{-\frac{\nu}{\lambda_{2}^j}x}}{\lambda_{2}^j}, \frac{\mathrm{e}^{-\frac{\nu}{\lambda_{3}^j}x}}{\lambda_{3}^j}, \frac{\mathrm{e}^{-\frac{\nu}{\lambda_{1}^j}x}}{\lambda_{1}^j} \right\}.
\end{align}
In what follows, we denote $c_i$ positive constants.
Combining Young's inequality with Lemma \ref{lem_bound_f}, we obtain:
\begin{align}
    &\int_0^L \left|\frac{2a}{\Lambda^-_j} \mathrm{e}^{\frac{\nu}{\Lambda^-_j} x} \beta(x,t)\mathbf{f}_{1j}(x)\mathbf{T}^+_j\mathbf{z}(x,t) \right|dx \notag \\
    &\leq \frac{a\mathrm{e}^{\frac{\nu}{\Lambda_j^-}}}{\underline{\Lambda}^-}M_0\left|s_2^0 - s_2^j\right|\left(  \int_0^L\beta^2(x,t) + 
    w_{1}^2(x,t) \right.\notag \\
    & +w_{2}^2(x,t) 
    +w_{3}^2(x,t)
    d x \Bigg) \notag \\
    &\leq c_1\left|s_2^0 - s_2^j\right|V(t),
\end{align}
where we have used the boundedness of the exponential term and the equivalence between the norm of the states $\vartheta$, $\mathbf{w}$, $\mathbf{z}$, and the Lyapunov functional $V_j$.  Consider now the term multiplied by $\mathbf{f}_{2j}(x)$. Using Young's inequality and Lemma \ref{lem_bound_f}, we get
\begin{align}
    &\int_0^L\left|\frac{2 a}{\Lambda^-_j} \mathrm{e}^{\frac{\nu}{\Lambda^-_j} x} \beta(x,t) \mathbf{f}_{2j}(x) \beta(0, t)\right| d x \notag \\
    &\leq \frac{2a\mathrm{e}^{\frac{\nu}{\Lambda_j^-}}}{\underline{\Lambda}^-} M_0|s_2^0 - s_2^j| \left( \int_0^L \beta(x,t)\beta(0,t)dx\right) \notag \\
    &\leq \frac{c_2}{k_1\varepsilon_0}\left|s_2^0 - s_2^j\right|V(t)+ c_2\varepsilon_0\left|s_2^0 - s_2^j\right| \beta^2(0,t). 
\end{align}
For the term  $\mathbf{f}_{3j}(x,\xi)$, we have:
\begin{align}
&\int_0^L \left|\frac{2a}{\Lambda^-_j} \mathrm{e}^{\frac{\nu}{\Lambda^-_j} x} \beta(x,t)\int_0^x\mathbf{f}_{3j}(x,\xi)\mathbf{T}^+_j\mathbf{z}(\xi,t) d\xi \right|dx \notag \\
&\leq c_3 \left|s_2^0-s_2^j\right|V(t).
\end{align}
Similarly, for the term $\mathbf{f}_{4j}(x,\xi)$, we get:
\begin{align}
&\int_0^L \left|\frac{2a}{\Lambda^-_j} \mathrm{e}^{\frac{\nu}{\Lambda^-_j} x} \beta(x,t)\int_0^x\mathbf{f}_{4j}(x,\xi)\mathbf{T}^-_j\mathbf{z}(\xi,t) d\xi \right|dx \nonumber \\
&\leq c_4 \left|s_2^0-s_2^j\right|V(t).
\end{align}
Therefore, we obtain
\begin{align}
&\frac{d V_j}{d \mathbf{z}}(\mathbf{z}) h_j(\mathbf{z})  \leq-\eta V_j(t)+M_1 \left|s_2^0-s_2^j\right| V(t) 
\notag \\
&+\left(c_2\left|\bar s-\underline s\right|\varepsilon_0+q_{1j}^2+q_{2j}^2+q_{3j}^2-a\right) \beta^2(0, t) \notag \\
&+\sum_{k=1}^3 (a \mathrm{e}^{\frac{\nu}{\Lambda^-_j}}\left(({R}_j)_k-({R}_0)_k\right)^2 -\mathrm{e}^{-\frac{\nu}{\lambda_{kj}}L})\alpha_k^2(L,t),
\end{align}
where 
\begin{align}
&\eta = \nu - \frac{2}{||\underline{\Lambda^+}|| k_1} (\max_{x\in [0,L]} ||\Sigma^{++}_0(x)||  \nonumber \\
    &+ (1 + \frac{1}{m_\vartheta}) \max_{x\in [0,L]} ||\Sigma^{+-}_0(x)||)-\frac{2\bar sc_2}{k_1\varepsilon_0}, \\
    &M_1=c_4+c_3+\frac{c_2}{k_1\varepsilon_0}+c_1.
\end{align}
The coefficients $\nu$, $\varepsilon_0$ and $a$ are chosen such that
\begin{align}
    &\eta >0, \quad c_2\left|\bar s-\underline s\right|\varepsilon_0+{q}_{1j}^2+{q}_{2j}^2+{q}_{3j}^2-a <0.
\end{align}
There exists a constant $C_0$ such that for all $1\leq j \leq r$, $V_j(\mathbf{z}) \leq C_0 V(\mathbf{z})$.
Thus, we get the following inequality:
\begin{align}
    &\frac{d V_j}{d \mathbf{z}}(\mathbf{z}) h_j(\mathbf{z})  \leq-\bar \eta V(t)+M_1 \left|s_2^0-s_2^j\right| V(t) \nonumber \\
    &+\sum_{k=1}^3 (a \mathrm{e}^{\frac{\nu}{\Lambda^-_j}}\left(({R}_j)_k-({R}_0)_k\right)^2 -\mathrm{e}^{-\frac{\nu}{\lambda_{kj}}L})\alpha_k^2(L,t),\label{ineq1}
\end{align}
where $\bar \eta=\eta C_0$.
Now, we calculate the second term of $L_j$. We have:
\begin{align}
    &\sum_{l=1}^r\left(V_l(\mathbf{z})-V_j(\mathbf{z})\right) \tau_{j l} =\sum_{l=1}^r  \tau_{j l}\notag \\
    &\left( \int_0^L \mathcal{K}_0(\mathbf{T}_l (\mathbf{z}(x,t)))^\mathsf{T} D_l(x) \mathcal{K}_0(\mathbf{T}_l (\mathbf{z}(x,t))) dx  \right.\notag\\
    &- \left.  \int_0^L \mathcal{K}_0(\mathbf{T}_j (\mathbf{z}(x,t)))^\mathsf{T} D_j(x) \mathcal{K}_0(\mathbf{T}_j (\mathbf{z}(x,t))) dx  \right).
\end{align}
For the integral term, we get:
\begin{align}
    &\int_0^L \mathcal{K}_0(\mathbf{T}_l (\mathbf{z}(x,t)))^\mathsf{T} D_l(x)  \mathcal{K}_0(\mathbf{T}_l (\mathbf{z}(x,t)))dx \notag\\
    &- \int_0^L \mathcal{K}_0(\mathbf{T}_j (\mathbf{z}(x,t)))^\mathsf{T} D_j(x)  \mathcal{K}_0(\mathbf{T}_j (\mathbf{z}(x,t))) dx \notag\\
    & = \int_0^L \mathcal{K}_0(\mathbf{T}_l (\mathbf{z}(x,t)))^\mathsf{T} (D_l(x) -D_j(x) ) \mathcal{K}_0(\mathbf{T}_l (\mathbf{z}(x,t)))\notag\\
    & + \mathcal{K}_0(\mathbf{T}_l(\mathbf{z}(x,t)))^\mathsf{T} D_j(x)  \mathcal{K}_0((\mathbf{T}_l - \mathbf{T}_j)(\mathbf{z}(x,t))) \notag\\
    &- \mathcal{K}_0((\mathbf{T}_j-\mathbf{T}_l) (\mathbf{z}(x,t)))^\mathsf{T} D_j(x)  \mathcal{K}_0(\mathbf{T}_j (\mathbf{z}(x,t))) dx.
\end{align}
Let us focus on the first term. We can apply the mean-value theorem to obtain
\begin{align}
    \int_0^L &\mathcal{K}_0(\mathbf{T}_j (\mathbf{z}(x,t)))^\mathsf{T} ( D_l(x)-D_j(x) ) \mathcal{K}_0(\mathbf{T}_j (\mathbf{z}(x,t))) dx \notag \\
    &\leq c_5 \left|s_2^{l}-s_2^j\right| V(t).
\end{align}
For the second term, we obtain
\begin{align}
   &\int_0^L \mathcal{K}_0(\mathbf{T}_l(\mathbf{z}(x,t)))^\mathsf{T} D_j(x) \mathcal{K}_0((\mathbf{T}_l - \mathbf{T}_j)(\mathbf{z}(x,t))) dx \notag\\
   & \leq \int_0^L||D_j(x)|| ||\mathbf{T}_l||  || \mathbf{T}_l -\mathbf{T}_j|| || \mathcal{K}_0||^2 ||\mathbf{z}||^2 dx\notag \\
   & \leq c_6 \left|s_2^{l}-s_2^j\right| V(t).
\end{align}
Similarly, we obtain for the third term
\begin{align}
    &- \int_0^L \mathcal{K}_0((\mathbf{T}_j-\mathbf{T}_l) (\mathbf{z}(x,t)))^\mathsf{T} D_j \mathcal{K}_0(\mathbf{T}_j (\mathbf{z}(x,t))) dx \notag\\
    & \leq c_7 \left|s_2^{l}-s_2^j\right| V(t).
\end{align}
Consequently, we get:
\begin{align}
    \sum_{l=1}^r\left(V_{l}(\vartheta)-V_j(\vartheta)\right) \tau_{jl} \leq  d_1 \sum_{l=1}^r \tau_{jl}\left|s_2^{l}-s_2^j\right| V(t), \label{ineq2}
\end{align}
where $d_1 = c_5 + c_6+ c_7$.
Using \eqref{ineq1} and \eqref{ineq2}
\begin{align}
    &L_j V(t) \leq -\bar \eta V(t)+M_1 \left|s_2^0-s_2^j\right| V(t) \notag\\
    &+ d_1 \sum_{l=1}^r  \tau_{j l}\left|s_2^l-s_2^j\right| V(t) \notag \\
    &+\sum_{k=1}^3 (a \mathrm{e}^{\frac{\nu}{\Lambda^-_j}}\left(({R}_j)_k-({R}_0)_k\right)^2 -\mathrm{e}^{-\frac{\nu}{\lambda_{kj}}L})\alpha_k^2(L,t).
\end{align}
We calculate the quantity $\Bar{L} = \sum_{j=1}^r P_{ij}(0,t) L_jV(t)$. Using the property of the expectation, we know that $\sum_{j=1}^r P_{i j}(0, t) |s_2^0-s_2^j|=\mathbb{E}\left(| s_2^0-s_2(t)|\right)$. Consequently, we obtain
\begin{align}
    &\Bar{L} = \sum_{j=1}^r \left[ P_{ij}(0,t) L_jV(t) 
    \leq \sum_{j=1}^r P_{i j}(0, t)\left(-\bar \eta V(t) \notag \right.\right.\\
    &+M_1 \left|s_2^0-s_2^j\right| V(t)+d_1 \sum_{l=1}^rP_{ij}(0,t)  \tau_{j l}\left|s_2^l-s_2^j\right| V(t) \notag \\
    &+P_{ij}(0,t)\sum_{k=1}^3 (a \mathrm{e}^{\frac{\nu}{\Lambda^-_j}}\left(({R}_j)_k-({R}_0)_k\right)^2 \notag\\
    &\left.-\mathrm{e}^{-\frac{\nu}{\lambda_{kj}}L})
    \alpha_k^2(L,t)\right], \notag\\
    &\leq -\bar \eta V(t)+ M_1 \mathbb{E}(|s_2^0 - s_2(t)|)V(t) \notag \\
    & + d_1 \sum_{j=1}^r P_{i j}(0, t) \sum_{l=1}^r  \tau_{jl}\left(\left| s_l-s_2^0\right|+\left|s_2^j-s_2^0\right|\right) V(t) \notag \\
    &+\sum_{k=1}^3 (d_2 \mathbb{E}(\left|s_2^0 - s_2(t)\right|)-\mathrm{e}^{-\frac{\nu}{\bar \lambda}L})\alpha_k^2(L,t),
\end{align}
where $\Bar{\lambda} = \max\{\lambda_{2}^j,\lambda_{3}^j,\lambda_{1}^j\}$, and we have used the fact that $\left|\left|(R_j)_k-(R_0)_k\right|\right|^2  \leq 2c_{10}\bar s\left|s_2^j -s_2^0\right|$. Consequently, we have:
\begin{align}
    &\Bar{L} \leq -V(t)\left(\bar \eta - (M_1+d_1 r \tau^\star) \mathbb{E}(\left|s_2^0-s_2(t)\right|) \right. \notag\\
    & \left.+ d_1 \sum_{j=1}^r\left|s_2^j - s_2^0\right|\left(\partial_t P_{i j}(0, t)+\mathfrak{c}_j P_{i j}(0, t)\right)\right) \notag \\
    &+\sum_{k=1}^3 (d_2 \mathbb{E}(\left|s_2^0 - s_2(t)\right|)-\mathrm{e}^{-\frac{\nu}{\bar \lambda}})\alpha_k^2(L,t),
\end{align}
This concludes the proof of Lemma \ref{lem_lyapunov_functional}. $\hfill\blacksquare$
\end{pf}
\subsection{Proof of Theorem \ref{main_thm}}
We now have all the tools to prove Theorem \ref{main_thm}. Notice first that if $\epsilon^\star$ is small enough (namely smaller than $\frac{e^{-\frac{\nu}{\bar \lambda}L}}{d_2}$) and if inequality~\eqref{eq_ineg_stocha} holds, then we have 
\begin{align}
&\sum_{j=1}^r P_{i j}(0, t) L_j V(t) \leq -V(t)\Bigg(\Bar{\eta}-d_1 \mathcal{Z}(t) \notag\\
&\left.-\left(M_1+d_1 r \tau^{\star}\right) \mathbb{E}\left(\left|s_2^0-s_2(t)\right|\right)\right).
\end{align}
We define the following function:
\begin{align}
    \phi(t) = \Bar{\eta }- d_1 \mathcal{Z}(t) - (M_1 + d_1 r \tau^\star)\mathbb{E}(|s_2^0-s_2(t)|).
\end{align}
And then, using the functional $\Psi(t)$:
\begin{align}
    \Psi(t) = \mathrm{e}^{\int_0^t \phi(y) dy} V(t).
\end{align}
With the definition of $\Psi(t)$, taking the expectation of the infinitesimal generator $L$ of $\Psi(t)$ , we get:
\begin{align}
    \mathbb{E}\left(\sum_{j=1}^r P_{ij}(0,t) L_{j}V(t)\right)\leq - \mathbb{E}\left(V(t)\phi(t)\right).
\end{align}
We know that $\mathbb{E}(\sum_{j=1}^r P_{ij}(0,t) L_{j}V(t)) =  \mathbb{E}(LV(t))$, thus
\begin{align}
    \mathbb{E}(LV(t))  \leq - \mathbb{E}(V(t)\phi(t)). \label{inequ}
\end{align}
Based on \eqref{inequ}, we get:
\begin{align}
    \mathbb{E}(L\Psi(t)) \leq 0.
\end{align}
Then applying the Dynkin's formula~\cite{dynkin2012theory},
\begin{align}
    \mathbb{E}(\Psi(t)) - \Psi(0) = \mathbb{E}\left(\int_0^t L\Psi(y)dy\right) \leq 0.
\end{align}
To calculate the $\mathbb{E}(\Psi(t))$, we write down the formulation of $\Psi(t)$:
\begin{align}
    &\mathbb{E}(\Psi(t)) =  \mathbb{E}\left( V(t) \mathrm{e}^{\int_0^t \phi(y)dy} \right) \notag \\
    & = \mathbb{E}\left( V(t) \mathrm{e}^{\int_0^t (\Bar{\eta} - d_1 \mathcal{Z}(y) - (M_1 + d_1 r \tau^*) \mathbb{E}(|s_2^0-s_2(y)|))dy} \right).
\end{align}
We already know that 
\begin{align}
    &\int_0^t \mathcal{Z}(y)dy =  \int_0^t \left(\sum_{j=1}^r\left|s_2^0 - s(y)\right|\left(\partial_y P_{i j}(0, y) \right.\right.\notag \\
    &\left. +\mathfrak{c}_j P_{i j}(0, y)\right) V(y)\Big) dy \notag\\
    &\leq \mathbb{E}(\left|s_2^0-s_2(t)\right|) + r\tau^\star \int_0^t \mathbb{E} (\left|s_2^0-s(y)\right|) dy,
\end{align}
where $\tau^\star$ is the largest value of the transition rate. Using this inequality, we get
\begin{align}
    &\mathbb{E}(\Psi(t)) \geq \mathbb{E}\left( V(t) \exp\left(\int_0^t (\Bar{\eta} - (M_1 + d_1 r \tau^\star) \right. \right. \notag\\
    & \mathbb{E}(\left|s_2^0-s(y)\right|))dy - d_1\mathbb{E}(\left|s_2^0-s(y)\right|) \notag \\
    & \left.\left.- d_1 r \tau^\star\int_0^t \mathbb{E} (\left|s_2^0-s(y)\right|) dy\right) \right), \notag \\
    &\geq \mathbb{E} \left( V(t) \mathrm{e}^{\left(-d_1 \epsilon^\star + \int_0^t (\Bar{\eta} - (M_1 + 2d_1 r \tau^\star) \epsilon^\star dy\right)} \right).
\end{align}
Then we take $\epsilon^\star$ as 
\begin{align}
    \epsilon^\star = \frac{\Bar{\eta}}{2(2d_1 r \tau^\star + M_1)},
\end{align}
thus we have 
\begin{align}
    \mathbb{E}(\Psi(t)) \geq \mathbb{E} \left( V(t) \mathrm{e}^{\left(-d_1 \epsilon^\star + \frac{\Bar{\eta}}{2} t\right)} \right).
\end{align}
From the before proof, we know $\mathbb{E}(\Psi(t)) \leq \Psi(0)$, such that
\begin{align}
    \mathbb{E}(V(t)) \leq  \mathrm{e}^{d_1 \epsilon^\star} \mathrm{e}^{-\zeta t} V(0),
\end{align}
where $\zeta = \frac{\Bar{\eta}}{2}$. The function $V(t)$ is equivalent to the $L^2$-norm of the system. This concludes the proof of Theorem \ref{main_thm}. $\hfill\blacksquare$

\section{Numerical Simulation}
In this section, we illustrate our results with simulations. The density and velocity evolution of HVs and AVs are discussed, and the probability of stochastic parameters and the states that are reached in the simulation process are presented. 
\subsection{Simulation configuration}
\begin{figure}
\centering
    \begin{subfigure}{0.23\textwidth}
        \includegraphics[width =\textwidth]{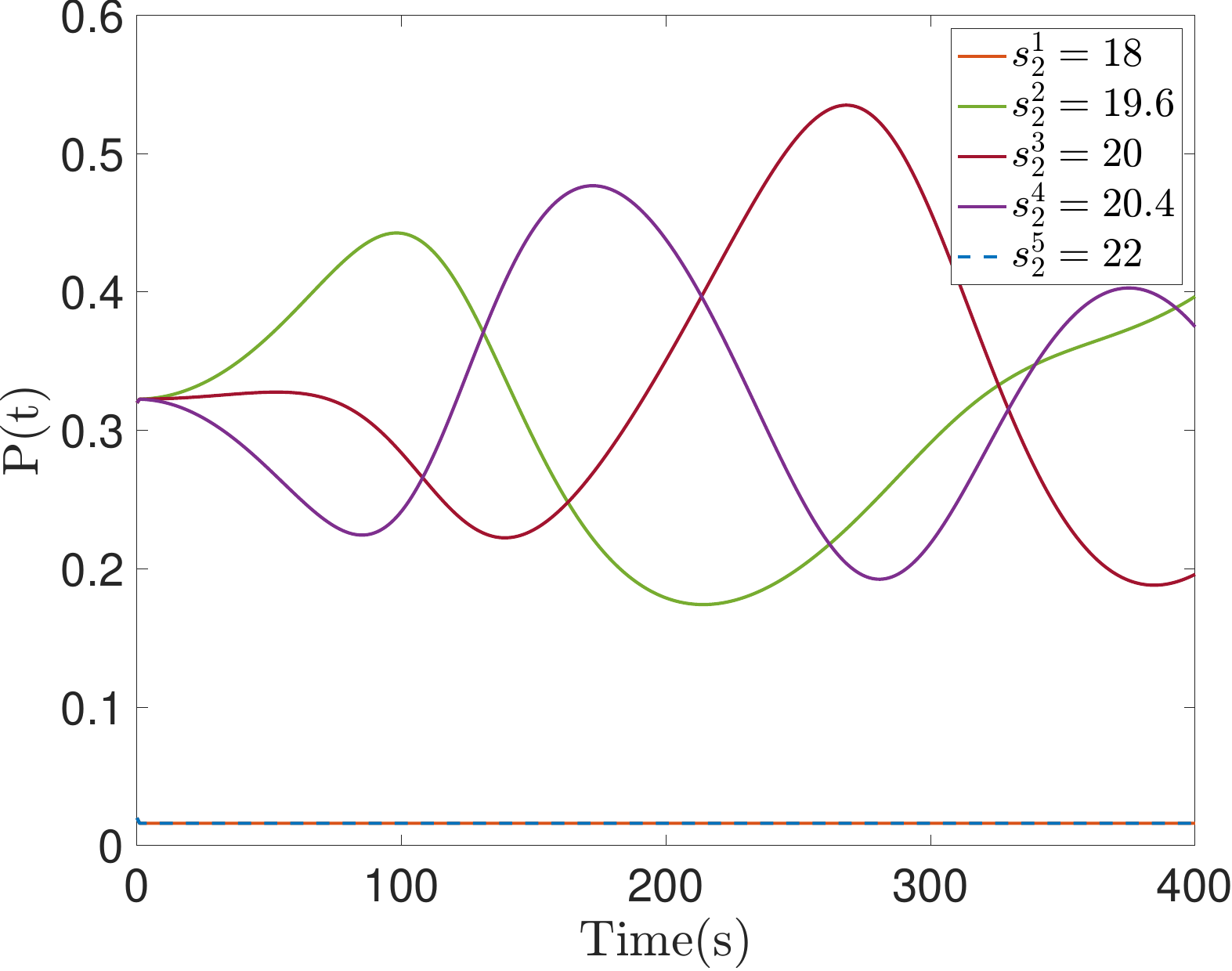}
        \caption{Probability of states}
        \label{Pro}
    \end{subfigure}
    \begin{subfigure}{0.23\textwidth}
        \includegraphics[width =\textwidth]{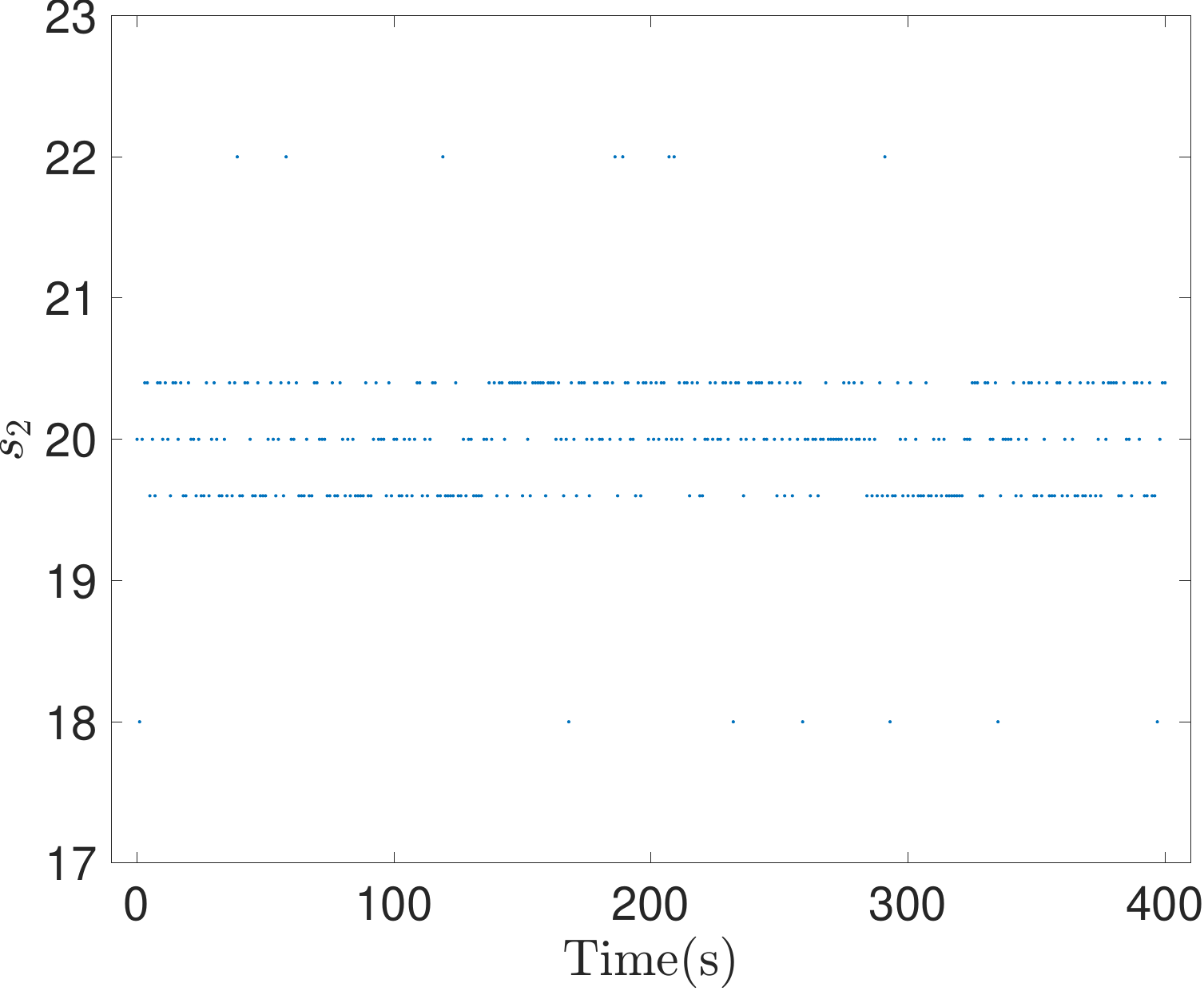}
        \caption{States reached}
        \label{states_reached}
    \end{subfigure}
    \caption{The probability of Markov states and the states the system reached}
    \label{trpro}
\end{figure}

In the implementation of the simulation, we set the nominal equilibrium density $\rho_1^* = 150 \text{veh/km}$, $\rho_2^*=75\text{veh/km}$. The corresponding nominal equilibrium velocity are $v_1^* = 29.16\text{km/h}$, $v_2^* = 13.32\text{km/h}$. The free-flow velocity for HVs and AVs are $V_1 = 80\text{km/h}$, $V_2 = 60\text{km/h}$ and the relaxation time for them are chosen as $\iota_1 = 30\text{s}$, $\iota_2 = 60\text{s}$. The pressure exponent value is selected as $\gamma_1 = 2.5$, $\gamma_2 = 2$. The nominal value of area occupancy are $a_1 = 10 \text{m}^2$, $a_2 = 40 \text{m}^2$, indicating the nominal spacing are $s_1 = 5 \text{m}$, $s_2 = 20 \text{m}$. The corresponding maximum area occupancy are $\overline{AO}_1 = 0.9$, $\overline{AO}_2 = 0.85$. We simulate on a $ L = 1000\text{m}$ long road whose width is $6 \text{m}$ with two lanes. The simulation time is $400\text{s}$. We consider the spacing of AVs $s_2$ as stochastic. Following the basic properties of a continuous Markov chain in the previous section, we choose five different values for $s_2$. The five values for spacing of AVs are $s_2^1=18, s_2^2=19.6, s_2^3=20, s_2^4=20.4, s_2^5=22)$ and the initial transition probabilities are chosen as $(0.02,0.32,0.32,0.32,0.02)$. 
The transition rates $\tau_{ij}$ are defined from~\cite{auriol2023mean}:
\begin{align}
    \tau_{i j}(t)= \begin{cases}0, \quad \text { if } i=j \\ 2 \mathfrak{r}, \quad \text { if } i \in\{1,5\} \\ 0.1 \mathfrak{r}, \quad  \text { if } i \in\{2,3,4\}, j \in\{1,5\} \\ \mathfrak{r}\left(1+2 \cos (\mathfrak{s}(i+5 j) t)^2\right), \text{others}  \end{cases}
\end{align}
where $\mathfrak{r}=10$ and $\mathfrak{s} = 0.001$. Solving the Kolmogorov forward equation, we get the probability of each state in the simulation process shown in Fig. \ref{Pro}, and the states reached in the process are shown in Fig. \ref{states_reached}. From the probability of each state, the system is more likely to jump to the values that are close to the nominal value. As shown in Fig. \ref{states_reached}, the system stays in $\{19.6, 20, 20.4\}$ more often in the simulation process. The sinusoidal initial conditions are applied to describe the stop-and-go traffic scenarios.
\subsection{Simulation results}
\begin{figure}
\centering
    \begin{subfigure}{0.4\textwidth}
        \includegraphics[width =\textwidth]{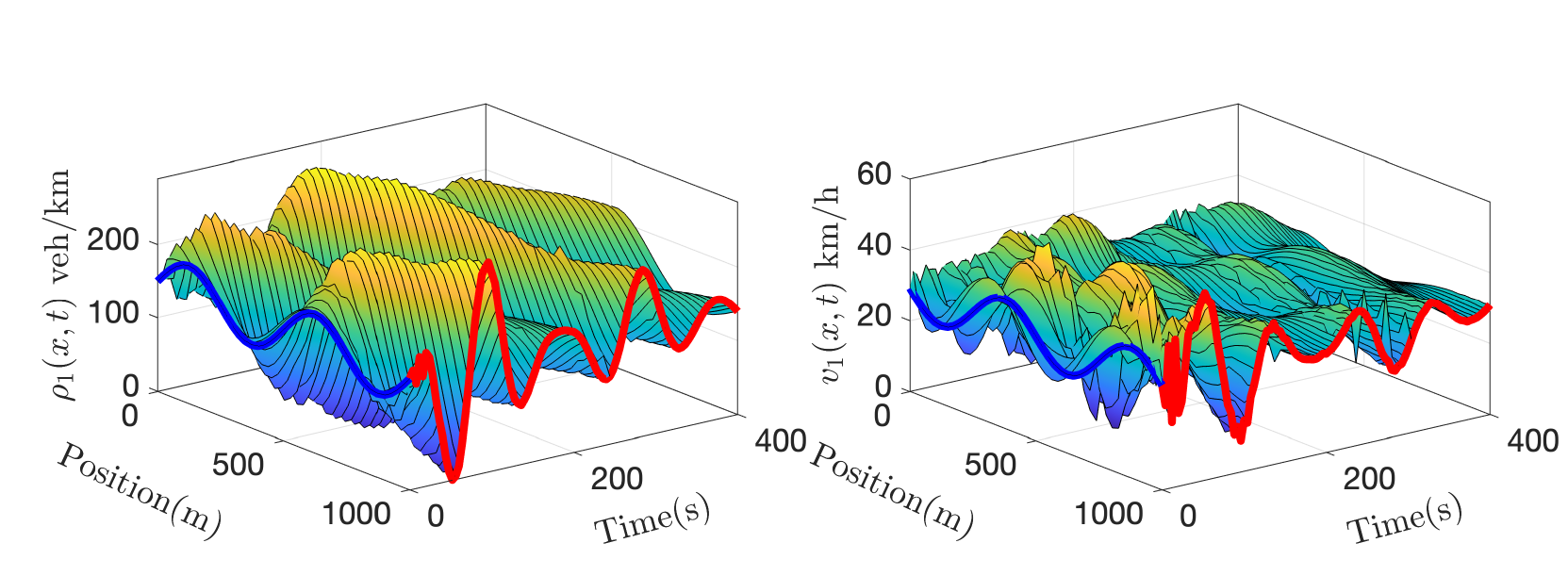}
        \caption{HVs}
        \label{ol_Markov_c1}
    \end{subfigure}
    \begin{subfigure}{0.4\textwidth}
        \includegraphics[width =\textwidth]{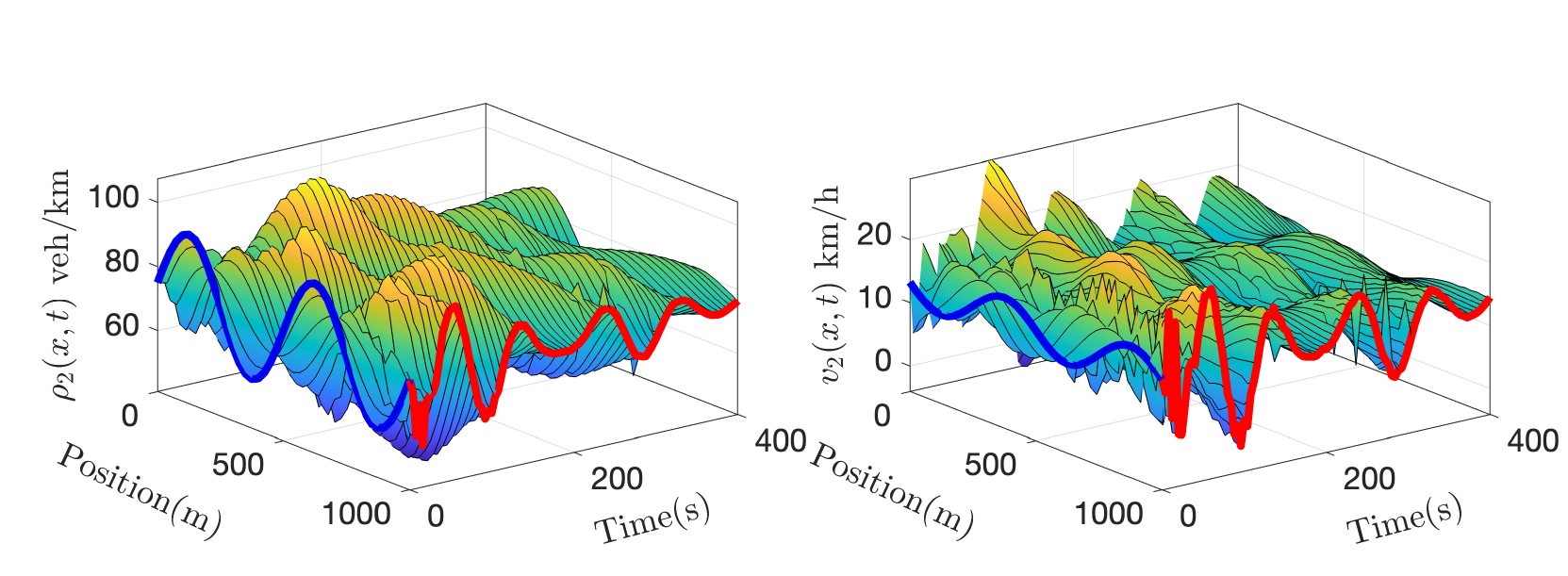}
        \caption{AVs}
        \label{ol_Markov_c2}
    \end{subfigure}
    \caption{The open-loop results of the stochastic system}
    \label{ol_Markov}
\end{figure}

\begin{figure}
\centering
    \begin{subfigure}{0.4\textwidth}
        \includegraphics[width =\textwidth]{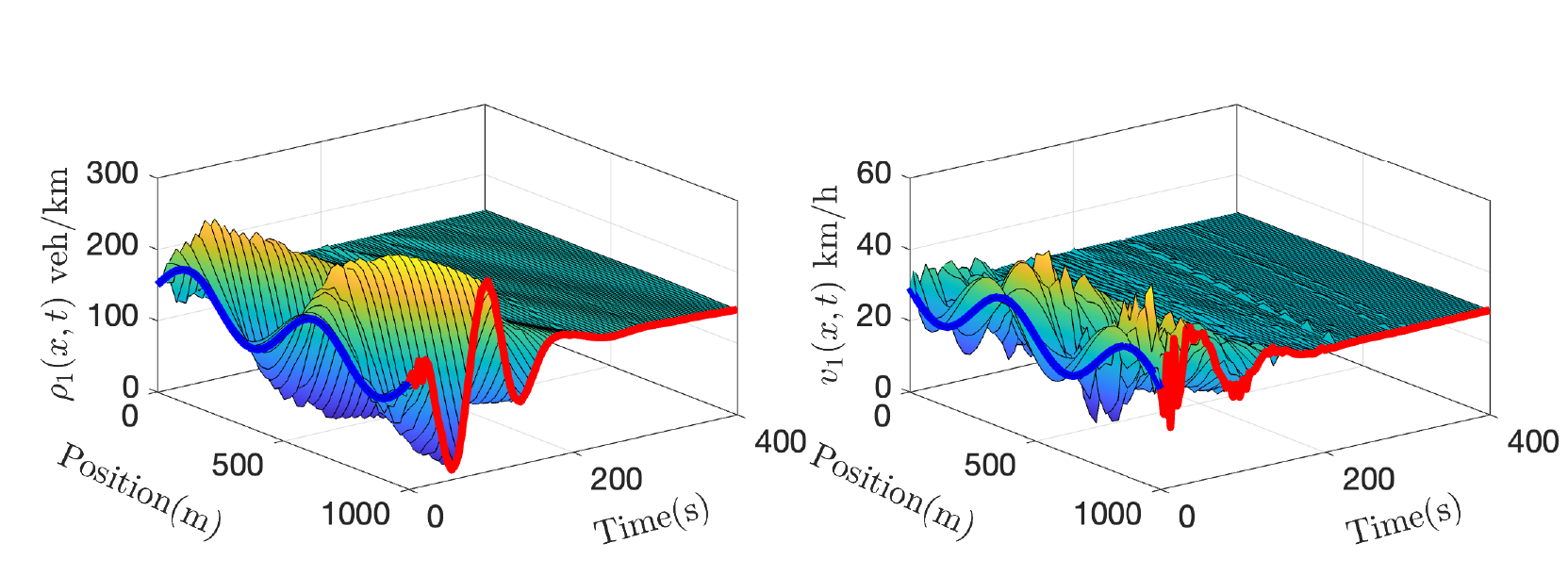}
        \caption{HVs}
        \label{cl_Markov_c1}
    \end{subfigure}
    \begin{subfigure}{0.4\textwidth}
        \includegraphics[width =\textwidth]{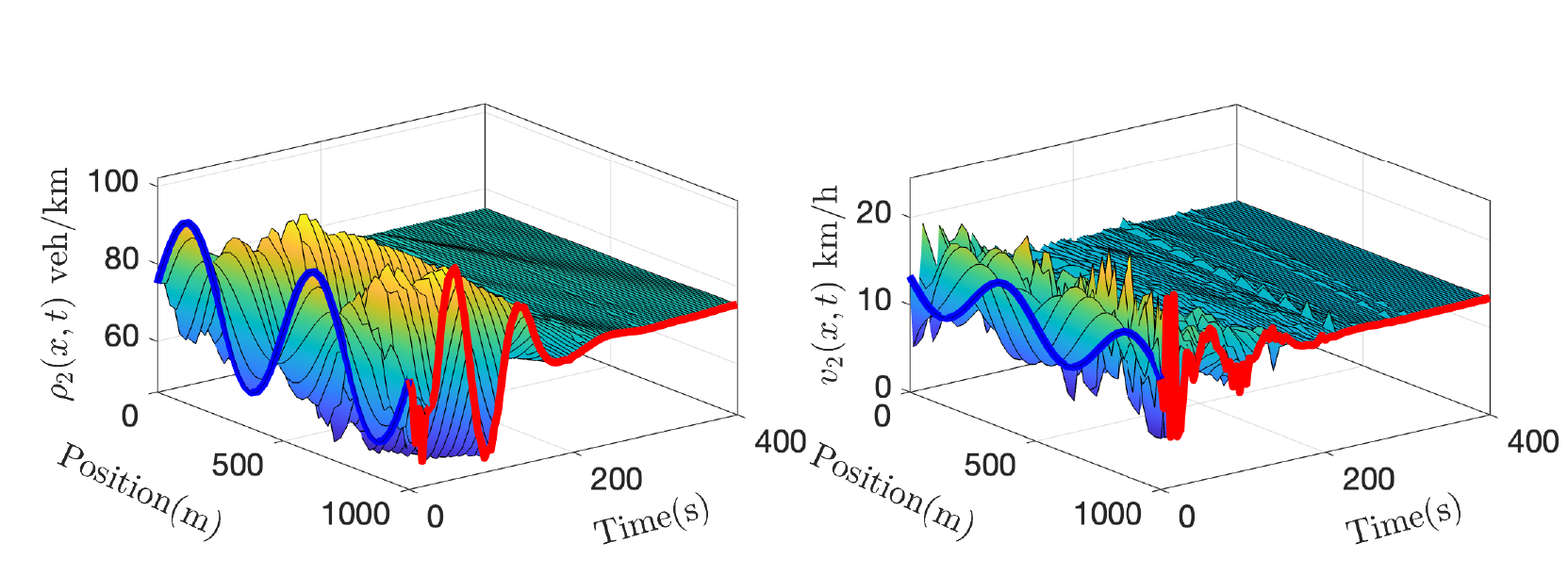}
        \caption{AVs}
        \label{cl_Markov_c2}
    \end{subfigure}
    \caption{The close-loop results of the stochastic system}
    \label{cl_Markov}
\end{figure}

\begin{figure}
\centering
    \begin{subfigure}{0.45\textwidth}
        \includegraphics[width =\textwidth]{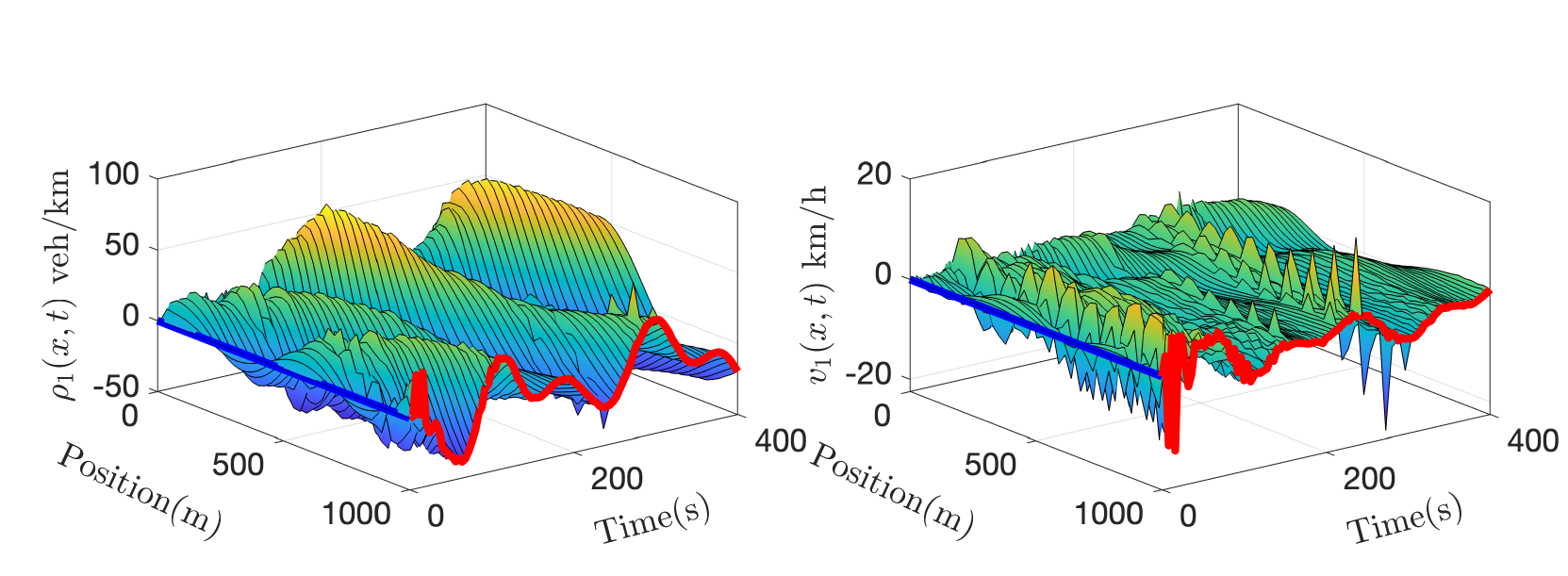}
        \caption{HVs}
        \label{ol_com_c1}
    \end{subfigure}
    \begin{subfigure}{0.45\textwidth}
        \includegraphics[width =\textwidth]{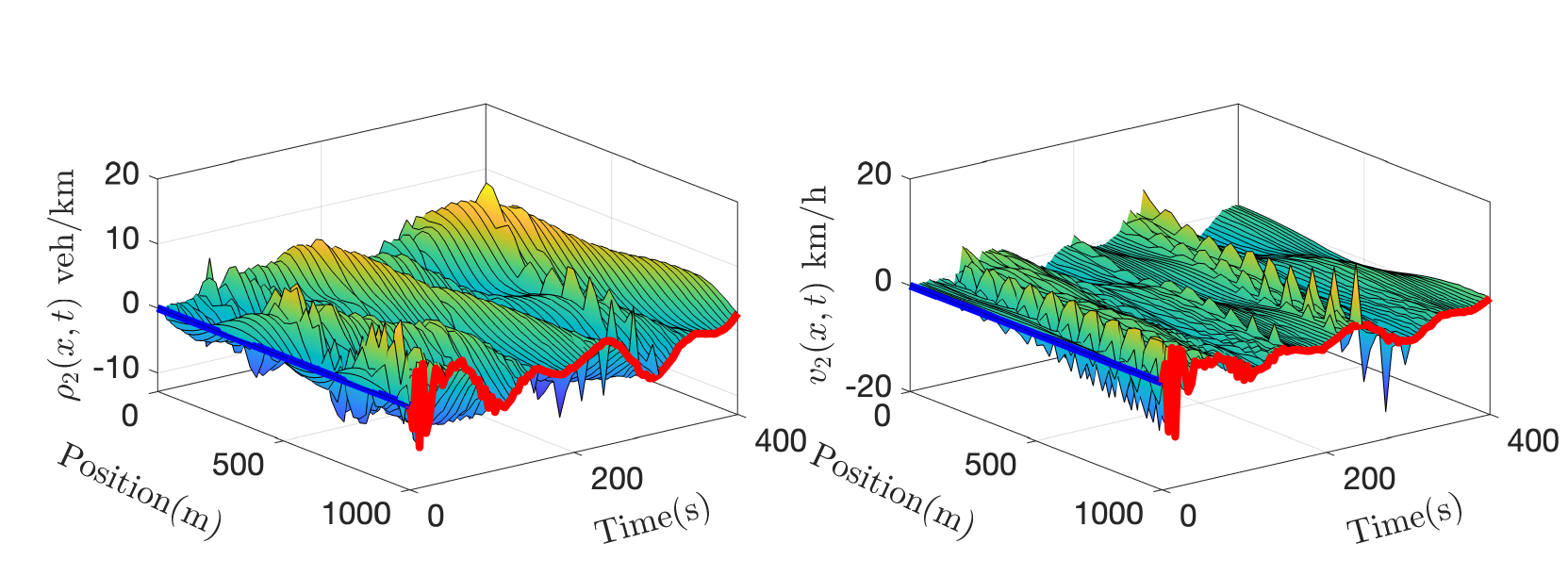}
        \caption{AVs}
        \label{ol_com_c2}
    \end{subfigure}
    \caption{The density and velocity error between open-loop stochastic system and nominal system}
    \label{ol_com}
\end{figure}

\begin{figure}
\centering
    \begin{subfigure}{0.45\textwidth}
        \includegraphics[width =\textwidth]{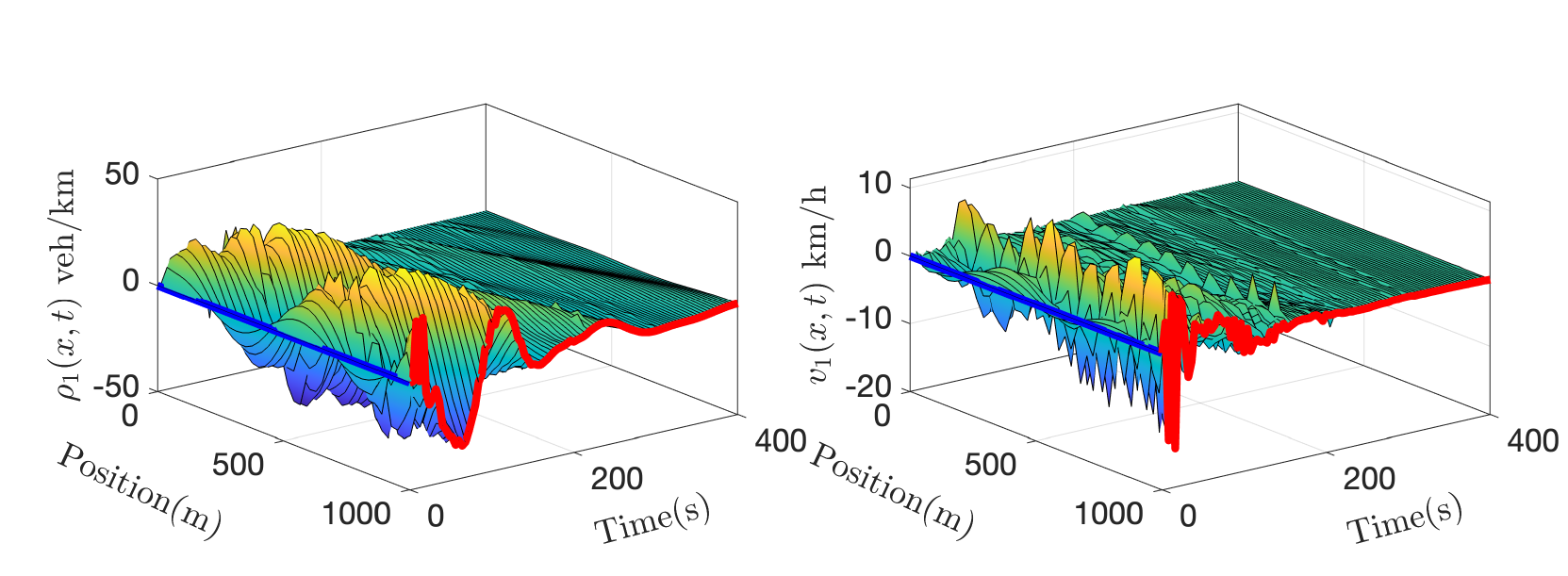}
        \caption{HVs}
        \label{cl_com_c1}
    \end{subfigure}
    \begin{subfigure}{0.45\textwidth}
        \includegraphics[width =\textwidth]{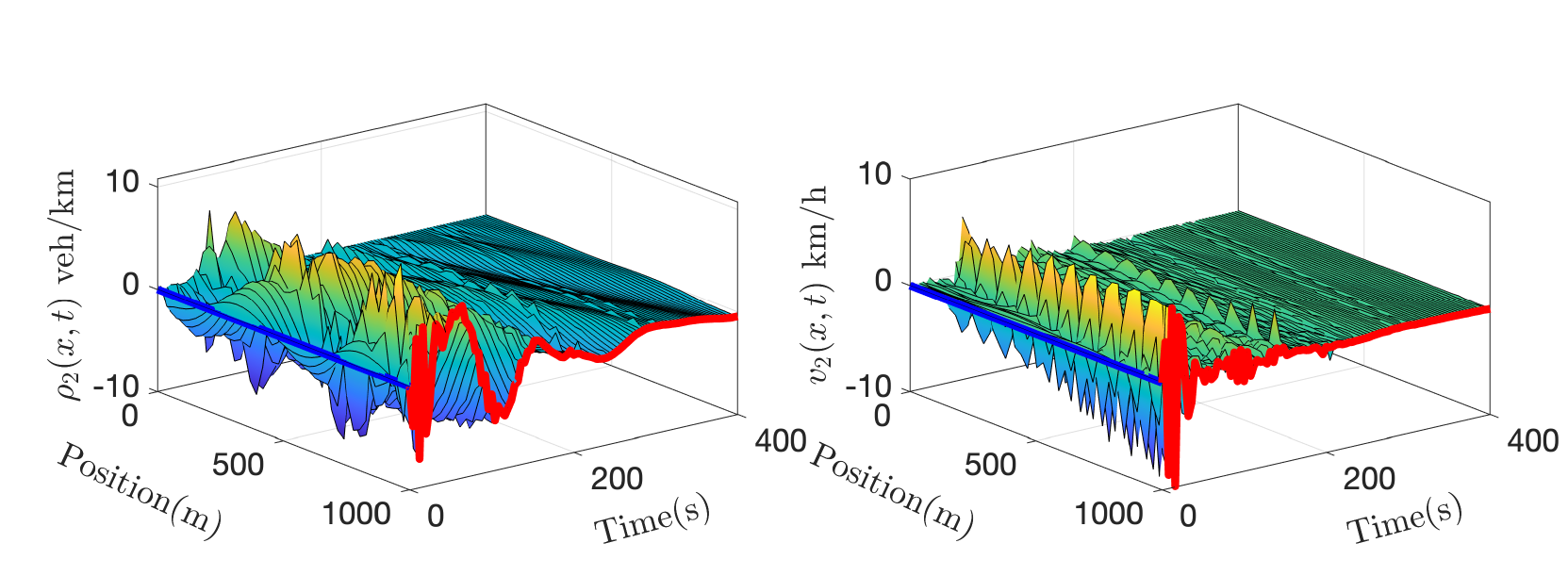}
        \caption{AVs}
        \label{cl_com_c2}
    \end{subfigure}
    \caption{The density and velocity error between closed-loop stochastic system and nominal system}
    \label{cl_com}
\end{figure}

\begin{figure}
    \centering
    \includegraphics[width = 0.45\textwidth]{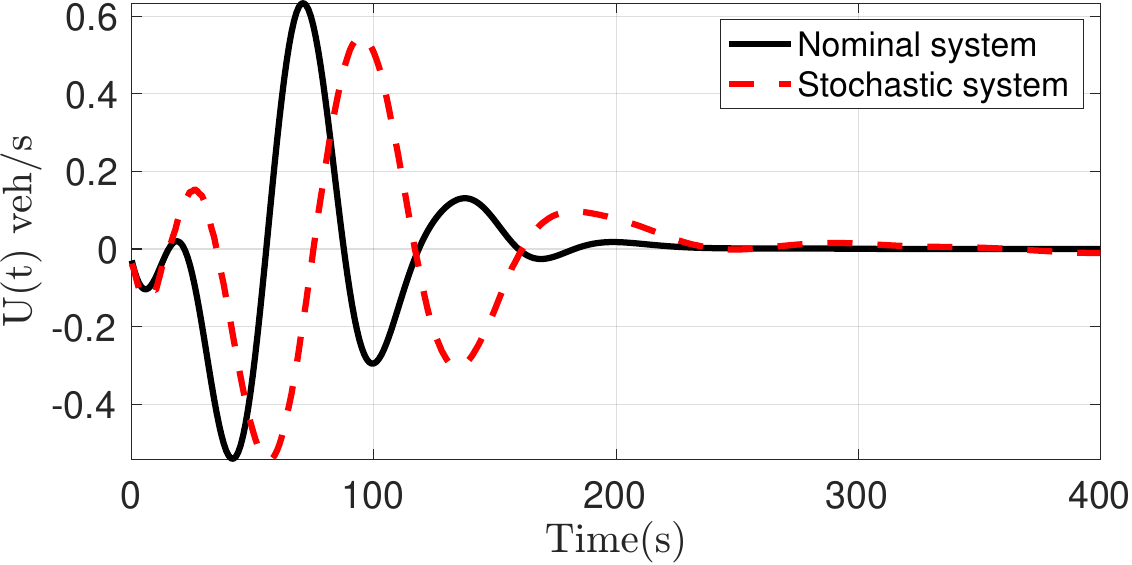}
    \caption{The control law for nominal system and stochastic system}
    \label{control_law_fig}
\end{figure}

The open-loop behavior of the system with stochastic parameters is shown in Fig. \ref{ol_Markov_c1},\ref{ol_Markov_c2}. The density and velocity of the mixed-autonomy traffic without control oscillate throughout the whole process. The closed-loop behavior of our stochastic system is shown in Fig. \ref{cl_Markov_c1},\ref{cl_Markov_c2}.  
The traffic states deviate from the original states due to the stochastic driving strategy of AVs. The nominal system of HVs and AVs reach their equilibrium points at a finite time with the proposed backstepping control law.  
We also provide a comparison of the nominal system and stochastic system for open-loop and closed-loop conditions. The error for open-loop results between the nominal and stochastic systems are shown in Fig. \ref{ol_com_c1},\ref{ol_com_c2}. The error for closed-loop results between the nominal system and stochastic system is shown in Fig. \ref{cl_com_c1},\ref{cl_com_c2}. It is observed that small deviations of density and velocity from the equilibrium still exist after the finite time stabilization is achieved. The maximum density and velocity error for open-loop HVs are $60 \text{veh/km}$,$20 \text{km/h}$ while the maximum density and velocity error of AVs are $10 \text{veh/km}$ and $18 \text{km/h}$. The maximum density and velocity error for the close-loop of AVs is $22 \text{veh/km}$, $15 \text{km/h}$. For HVs, the maximum density and velocity error are $10 \text{veh/km}$, $10\text{km/h}$. The evolution of control input $U(t)$ is shown in Fig. \ref{control_law_fig}.

In congested mixed-autonomy traffic conditions where AVs employ a more conservative spacing strategy in comparison to HVs, the impact areas of AVs are larger than HVs. The proposed control strategy can ensure the traffic system is in stable condition if the driving strategy of AVs remains relatively consistent with their nominal driving strategy. In the presence of larger gaps between AVs, HVs tend to exhibit gap-filling behavior, colloquially referred to as the ``creeping effect'', wherein HVs surpass AVs within such congested regimes.

\section{Conclusion}
In this paper, a mixed-autonomy traffic PDE model with Markov jumping parameters is developed to characterize the evolution of traffic density and velocity on the road. The spacing of AVs is considered as a Markov process, rendering the system stochastic in nature. The exponential stability of the nominal system is demonstrated through the employment of a full-state feedback control law using the backstepping method. Subsequently, Lyapunov analysis is utilized to prove the mean-square exponential stability of the stochastic system under the nominal control law. Finally, a numerical simulation is conducted to assess the performance of the control approach for the stochastic system. It is of the author's future interest to develop a boundary observer and output-feedback boundary controllers to stabilize the stochastic system.

\section{Acknowledgments}
This work was supported in part by the
National Natural Science Foundation of China under Grant 62203131, 
in part by the Guangzhou Municipal Natural Science Foundation of under Grant 2023A04J01850. 

\bibliography{reference}           

\begin{thebibliography}{10}

\bibitem{albaba2019stochastic}
M.~Albaba, Y.~Yildiz, N.~Li, I.~Kolmanovsky, and A.~Girard.
\newblock Stochastic driver modeling and validation with traffic data.
\newblock In {\em 2019 American Control Conference (ACC)}, pages 4198--4203. IEEE, 2019.

\bibitem{amin2011exponential}
S.~Amin, F.~M. Hante, and A.~M. Bayen.
\newblock Exponential stability of switched linear hyperbolic initial-boundary value problems.
\newblock {\em IEEE Transactions on Automatic Control}, 57(2):291--301, 2011.

\bibitem{anfinsen2018adaptive}
H.~Anfinsen and O.~M. Aamo.
\newblock Adaptive control of linear 2$\times$ 2 hyperbolic systems.
\newblock {\em Automatica}, 87:69--82, 2018.

\bibitem{auriol2016minimum}
J.~Auriol and F.~Di~Meglio.
\newblock Minimum time control of heterodirectional linear coupled hyperbolic pdes.
\newblock {\em Automatica}, 71:300--307, 2016.

\bibitem{auriol2020robust}
J.~Auriol and F.~Di~Meglio.
\newblock Robust output feedback stabilization for two heterodirectional linear coupled hyperbolic {PDE}s.
\newblock {\em Automatica}, 115:108896, 2020.

\bibitem{auriol2023mean}
J.~Auriol, M.~Pereira, and B.~Kulcsar.
\newblock Mean-square exponential stabilization of coupled hyperbolic systems with random parameters.
\newblock {\em IFAC-PapersOnLine}, pages 8823--8828, 2023.

\bibitem{aw2000resurrection}
A.~Aw and M.~Rascle.
\newblock Resurrection of "second order" models of traffic flow.
\newblock {\em SIAM journal on applied mathematics}, 60(3):916--938, 2000.

\bibitem{bastin2016stability}
G.~Bastin and J.-M. Coron.
\newblock {\em Stability and boundary stabilization of 1-d hyperbolic systems}, volume~88.
\newblock Springer, 2016.

\bibitem{bekiaris2020pde}
N.~Bekiaris-Liberis and A.~I. Delis.
\newblock {PDE}-based feedback control of freeway traffic flow via time-gap manipulation of acc-equipped vehicles.
\newblock {\em IEEE Transactions on Control Systems Technology}, 29(1):461--469, 2020.

\bibitem{bolzern2006almost}
P.~Bolzern, P.~Colaneri, and G.~De~Nicolao.
\newblock On almost sure stability of continuous-time markov jump linear systems.
\newblock {\em Automatica}, 42(6):983--988, 2006.

\bibitem{burkhardt2021stop}
M.~Burkhardt, H.~Yu, and M.~Krstic.
\newblock Stop-and-go suppression in two-class congested traffic.
\newblock {\em Automatica}, 125:109381, 2021.

\bibitem{calvert2017will}
S.~Calvert, W.~Schakel, and J.~Van~Lint.
\newblock Will automated vehicles negatively impact traffic flow?
\newblock {\em Journal of advanced transportation}, 2017, 2017.

\bibitem{coron2021boundary}
J.-M. Coron, L.~Hu, G.~Olive, and P.~Shang.
\newblock Boundary stabilization in finite time of one-dimensional linear hyperbolic balance laws with coefficients depending on time and space.
\newblock {\em Journal of Differential Equations}, 271:1109--1170, 2021.

\bibitem{do2012continuous}
O.~L. do~Valle~Costa, M.~D. Fragoso, and M.~G. Todorov.
\newblock {\em Continuous-time Markov jump linear systems}.
\newblock Springer Science \& Business Media, 2012.

\bibitem{dynkin2012theory}
E.~B. Dynkin.
\newblock {\em Theory of Markov processes}.
\newblock Courier Corporation, 2012.

\bibitem{espitia2022traffic}
N.~Espitia, J.~Auriol, H.~Yu, and M.~Krstic.
\newblock Traffic flow control on cascaded roads by event-triggered output feedback.
\newblock {\em International Journal of Robust and Nonlinear Control}, 32(10):5919--5949, 2022.

\bibitem{hoyland2009system}
A.~Hoyland and M.~Rausand.
\newblock {\em System reliability theory: models and statistical methods}.
\newblock John Wiley \& Sons, 2009.

\bibitem{hu2015control}
L.~Hu, F.~Di~Meglio, R.~Vazquez, and M.~Krstic.
\newblock Control of homodirectional and general heterodirectional linear coupled hyperbolic {PDE}s.
\newblock {\em IEEE Transactions on Automatic Control}, 61(11):3301--3314, 2016.

\bibitem{jin2020analysis}
L.~Jin, M.~{\v{C}}i{\v{c}}i{\'c}, K.~H. Johansson, and S.~Amin.
\newblock Analysis and design of vehicle platooning operations on mixed-traffic highways.
\newblock {\em IEEE Transactions on Automatic Control}, 66(10):4715--4730, 2020.

\bibitem{kesting2008adaptive}
A.~Kesting, M.~Treiber, M.~Sch{\"o}nhof, and D.~Helbing.
\newblock Adaptive cruise control design for active congestion avoidance.
\newblock {\em Transportation Research Part C: Emerging Technologies}, 16(6):668--683, 2008.

\bibitem{kolmanovsky2001mean}
I.~Kolmanovsky and T.~L. Maizenberg.
\newblock Mean-square stability of nonlinear systems with time-varying, random delay.
\newblock {\em Stochastic analysis and Applications}, 19(2):279--293, 2001.

\bibitem{krstic2008boundary}
M.~Krstic and A.~Smyshlyaev.
\newblock {\em Boundary control of PDEs: A course on backstepping designs}.
\newblock SIAM, 2008.

\bibitem{lamare2015switching}
P.-O. Lamare, A.~Girard, and C.~Prieur.
\newblock Switching rules for stabilization of linear systems of conservation laws.
\newblock {\em SIAM Journal on Control and Optimization}, 53(3):1599--1624, 2015.

\bibitem{li2022trade}
X.~Li.
\newblock Trade-off between safety, mobility and stability in automated vehicle following control: An analytical method.
\newblock {\em Transportation research part B: methodological}, 166:1--18, 2022.

\bibitem{lighthill1955kinematic}
M.~J. Lighthill and G.~B. Whitham.
\newblock On kinematic waves ii. a theory of traffic flow on long crowded roads.
\newblock {\em Proceedings of the Royal Society of London. Series A. Mathematical and Physical Sciences}, 229(1178):317--345, 1955.

\bibitem{mohan2017heterogeneous}
R.~Mohan and G.~Ramadurai.
\newblock Heterogeneous traffic flow modelling using second-order macroscopic continuum model.
\newblock {\em Physics Letters A}, 381(3):115--123, 2017.

\bibitem{prieur2014stability}
C.~Prieur, A.~Girard, and E.~Witrant.
\newblock Stability of switched linear hyperbolic systems by lyapunov techniques.
\newblock {\em IEEE Transactions on Automatic control}, 59(8):2196--2202, 2014.

\bibitem{qi2022delay}
J.~Qi, S.~Mo, and M.~Krstic.
\newblock Delay-compensated distributed pde control of traffic with connected/automated vehicles.
\newblock {\em IEEE Transactions on Automatic Control}, 2022.

\bibitem{rausand2003system}
M.~Rausand and A.~Hoyland.
\newblock {\em System reliability theory: models, statistical methods, and applications}, volume 396.
\newblock John Wiley \& Sons, 2003.

\bibitem{redaud2024domain}
J.~Redaud, J.~Auriol, and Y.~Le~Gorrec.
\newblock In domain dissipation assignment of boundary controlled {P}ort-{H}amiltonian systems using backstepping.
\newblock {\em Systems \& Control Letters}, 185:105722, 2024.

\bibitem{richards1956shock}
P.~I. Richards.
\newblock Shock waves on the highway.
\newblock {\em Operations research}, 4(1):42--51, 1956.

\bibitem{ross2014introduction}
S.~M. Ross.
\newblock {\em Introduction to probability models}.
\newblock Academic press, 2014.

\bibitem{vazquez2011backstepping}
R.~Vazquez, M.~Krstic, and J.-M. Coron.
\newblock Backstepping boundary stabilization and state estimation of a 2$\times$ 2 linear hyperbolic system.
\newblock In {\em 2011 50th IEEE conference on decision and control and european control conference}, pages 4937--4942. IEEE, 2011.

\bibitem{wang2012stochastically}
J.-W. Wang, H.-N. Wu, and H.-X. Li.
\newblock Stochastically exponential stability and stabilization of uncertain linear hyperbolic pde systems with markov jumping parameters.
\newblock {\em Automatica}, 48(3):569--576, 2012.

\bibitem{yoshida1960lectures}
K.~Yoshida.
\newblock {\em Lectures on differential and integral equations}, volume~10.
\newblock Interscience Publishers, 1960.

\bibitem{yu2020stability}
H.~Yu, S.~Amin, and M.~Krstic.
\newblock Stability analysis of mixed-autonomy traffic with cav platoons using two-class {A}w-{R}ascle model.
\newblock In {\em 2020 59th IEEE Conference on Decision and Control (CDC)}, pages 5659--5664. IEEE, 2020.

\bibitem{yu2019traffic}
H.~Yu and M.~Krstic.
\newblock Traffic congestion control for {A}w--{R}ascle--{Z}hang model.
\newblock {\em Automatica}, 100:38--51, 2019.

\bibitem{yu2022traffic}
H.~Yu and M.~Krstic.
\newblock {\em Traffic Congestion Control by PDE Backstepping}.
\newblock Springer, 2022.

\bibitem{zhang2002non}
H.~M. Zhang.
\newblock A non-equilibrium traffic model devoid of gas-like behavior.
\newblock {\em Transportation Research Part B: Methodological}, 36(3):275--290, 2002.

\bibitem{zhang2017stochastic}
L.~Zhang and C.~Prieur.
\newblock Stochastic stability of markov jump hyperbolic systems with application to traffic flow control.
\newblock {\em Automatica}, 86:29--37, 2017.

\bibitem{zhang2019pi}
L.~Zhang, C.~Prieur, and J.~Qiao.
\newblock {PI} boundary control of linear hyperbolic balance laws with stabilization of arz traffic flow models.
\newblock {\em Systems \& Control Letters}, 123:85--91, 2019.

\bibitem{zhang2006hyperbolicity}
P.~Zhang, R.-X. Liu, S.~Wong, and S.-Q. Dai.
\newblock Hyperbolicity and kinematic waves of a class of multi-population partial differential equations.
\newblock {\em European Journal of Applied Mathematics}, 17(2):171--200, 2006.

\end{thebibliography}

\end{document}